\documentclass[12pt]{article}
\usepackage{amsmath,amsthm,amsfonts,latexsym,amscd,eucal}

\newcommand{\Reali}{{\mathbb{R}}}                                
\newcommand{\Complessi}{{\mathbb{C}}}                                 
                                 
\newcommand{\Interi}{{\mathbb{Z}}}                                 
                                 
\newcommand{\cO}{{\cal O}}                                      
\newcommand{\cC}{{\cal C}}                                      
\newcommand{\cA}{{\cal A}}                                      
\newcommand{\cB}{{\cal B}} 
                                      
\newcommand{\cH}{{\cal H}}

\newcommand{\cI}{{\cal I}}
\newcommand{\cN}{{\cal N}}
\newcommand{\cL}{{\cal L}}
\newcommand{\cM}{{\cal M}}

\newcommand{\cF}{{\cal F}}
\newcommand{\cR}{{\cal R}}

\newcommand{\cK}{{\cal K}}
\newcommand{\cT}{{\cal T}}
\newcommand{\cV}{{\cal V}}

\newcommand{\str}{\textnormal{Tr}_s}
\def\qed{{\unskip\nobreak\hfil\penalty50
\hskip2em\hbox{}\nobreak\hfil$\Box$
\parfillskip=0pt \finalhyphendemerits=0\par}\medskip}
\def\a{\alpha}
\def\b{\beta}
\def\d{\delta}
\def\D{\Delta}

\def\f{\varphi}
\def\g{\gamma}
\def\G{\Gamma}

\def\k{\kappa}

\def\l{\lambda}
\def\L{\Lambda}
\def\m{\mu}

\def\o{\omega}

\def\r{\rho}
\def\s{\sigma}

\def\t{\tau}

\def\x{\xi}

\def\gF{\mathfrak F}
\def\gA{\mathfrak A}
\def\gB{\mathfrak B}
\def\gC{\mathfrak C}

\def\Ad{\textnormal{Ad}}
\def\Ppo{{\cal P}_+^\uparrow}           

\def\ac{\underset{t\to i}{\textnormal {anal.cont.\,}}}
\def\aci{\underset{t\to -i}{\textnormal {anal.cont.\,}}}
\def\acb{\underset{t\to i\beta}{\textnormal {anal.cont.\,}}}
\newtheorem{Thm}{Theorem}[section]
\newtheorem{Cor}[Thm]{Corollary}
\newtheorem{Prop}[Thm]{Proposition}
\newtheorem{Lemma}[Thm]{Lemma}
\theoremstyle{definition}

\theoremstyle{remark}

\title{\huge  Notes~for~a~Quantum~Index~Theorem}
\author{Roberto Longo\\
{}\\
Dipartimento di Matematica, Universit\`a di Roma ``Tor
Vergata''\\ Via della Ricerca Scientifica, I--00133 Roma, Italy.\\
E-mail: longo@mat.uniroma2.it}
\footnotetext{Supported in 
part by MURST and GNAFA--INDAM.}
\date{March 6, 2000}
\begin{document}
\maketitle
\markboth{Quantum index theorem}{Roberto Longo}

\begin{abstract}
We view DHR superselection sectors with finite statistics 
as Quantum Field Theory analogs of 
elliptic operators where KMS functionals play the role of the 
trace composed with the heat kernel regularization. 
We extend our local holomorphic dimension formula 
and prove an analogue of the index 
theorem in the Quantum Field Theory context. The analytic index 
is the Jones index, more precisely the minimal dimension, and,
on a 4-dimensional spacetime, the DHR theorem gives the integrality 
of the index. We introduce the notion of holomorphic dimension; 
the geometric dimension is then
defined as the part of the holomorphic dimension 
which is symmetric under charge conjugation. 
We apply the AHKT theory
of chemical potential and we extend it to 
the low dimensional case, by using conformal field theory. 
Concerning Quantum Field Theory on curved spacetime,
the geometry of the manifold
enters in the expression for the dimension. If a quantum black hole 
is described by a spacetime with bifurcate Killing horizon and 
sectors are localizable on the horizon, the logarithm of the 
holomorphic dimension is proportional to the incremental free
energy, due to the addition of the charge, and to the inverse 
temperature, hence to the surface gravity in the Hartle-Hawking KMS 
state. For this analysis we consider a conformal net obtained by
restricting the field to the horizon (``holography''). Compared
with our previous work on Rindler spacetime, 
this result differs inasmuch as it concerns true black hole 
spacetimes, like the Schwarzschild-Kruskal manifold, and pertains to 
the entropy of the black hole itself, rather than of the outside 
system. An outlook concerns a possible relation with supersymmetry 
and noncommutative geometry.
\end{abstract}

\newpage
\tableofcontents 
\setcounter{section}{-1}
\section{Introduction.}
\label{sec:intro}
These notes are a natural outgrowth of our previous work on a local 
holomorphic formula for the dimension of a superselection sector
\cite{L4} and were motivated by the purpose to give a geometrical 
picture to aspects of local quantum physics related to the 
superselection structure. They may be read 
from different points of view, in particular, guided by a 
similarity of the statistical dimension with the Fredholm index
and a possible index theorem, already suggested in \cite{D}, 
we shall regard the DHR localized endomorphisms \cite{DHR} 
as quantum analogs of elliptic differential operators.

{\it A heuristic preamble.} 
We begin to give a heuristic, but elementary, motivation for our 
dimension formula, postponing for the moment the specification of the
underlying structure. Let the selfadjoint operator 
$H_0$ be a reference Hamiltonian for a 
Quantum Statistical Mechanics system, generating the evolution 
on an operator algebra $\gA$. $H_0$ may be thought to correspond
to the Laplacian $\D$ on a compact Riemann manifold.
We write any other Hamiltonian $H$ as
a perturbation $H=H_0 +P$ of $H_0$. We assume a supersymmetric 
structure, essentially the existence of a Dirac operator $D$ 
(an odd square root of $H$) that implements an odd derivation of $\gA$.

The McKean--Singer lemma then shows that
\begin{equation}
\str(e^{-\b H})= \text{Fredholm index of $D$}
		\label{MS}
\end{equation}
for any $\b>0$, where $\str$ denotes the super-trace.
In particular  $\str(e^{-\b H})$ (the Witten index) is an integer.

Now $\o^{(\b)} = 
{\str(e^{-\b H_{0}}\,\cdot)}/{\str(e^{-\b H_{0}})}$ is the 
normalized super-Gibbs functional at inverse temperature $\b$ 
for the dynamics generated by $H_0$ 
and if we consider the unitary cocycle
\begin{equation}
	u_P(t)\equiv e^{itH}e^{-itH_0}
	\label{coc}
\end{equation}
relating the two evolutions, and that belongs to $\gA$ if
$P\in\gA$, we may write the following formula
in terms of $\o^{(\b)}$ and  $u_P$:
\begin{equation}\label{ac1}
\underset{t\to i\b}{\textnormal {anal.cont.\,}}\o^{(\b)}(u_P(t))=
\frac{\str(e^{-\b H})}{\str(e^{-\b H_0})} \ ,
\end{equation}
provided the latter makes sense. We will thus 
regard the above expression as  a {\it multiplicative} 
relative index between $H_0$ and $H$.\footnote{As a counterpart,
an {\it additive} relative 
index generalizing $\str(e^{-\b H} - e^{-\b H_0})$ 
appeared  in \cite{BMS,J}.}

As shown in \cite{JLW,JL} this index is invariant under deformations,
in particular $\textnormal{Index}(u_P) = \textnormal{Index}(u_I)$ 
if $P$ is an odd element
of $\gA$ in the domain of the superderivation and may be obtained by
evaluating at the identity the JLO Chern character 
associated with $\o^{(\b)}$.

At infinite volume however, the integrality of the index is not evident.

In the case of infinite volume systems, likewise for the Laplacian on 
a non-compact manifolds, the Hamiltonian has not any longer discrete 
spectrum and may exist only as a derivation. What survives after 
the thermodynamical limit is 
the time evolution, a one-parameter automorphism group $\a$ of a 
C$^*$-algebra $\gA$. The equilibrium states $\o^{(\b)}$ 
are  characterized by 
the KMS condition \cite{HHW}. However the cocycle (\ref{coc}) may well 
exist and belong to $\gA$, so that formula (\ref{ac1}) may be generalized to
define {\it the index of a cocycle}, if the analytic 
continuation exists.

In this paper we shall consider the case of 
Relativistic Quantum Statistical Mechanics, by which we mean the 
consideration of KMS 
functionals for the time evolution in Quantum Field Theory.
The basic relativistic property, locality or finite 
propagation of speed of light, will select the appropriate class of 
cocycles and will imply the integrality of the index as will be 
explained. For reasons that will partly 
be clarified, the supersymmetric structure will not play a direct 
role in our formulae for the dimension.

{\it Ingredients for QFT analysis}.
Let us discuss now some fundamental aspects of Physics and Analysis. 
Quantum Field Theory can be considered at the same time 
as a generalization of two Physical theories of very 
different nature: Classical (Lagrangian) Field Theory and Quantum Mechanics. 
Both of them extend Classical Mechanics, but point in apparently 
divergent directions. In the first case one goes from finitely many to 
infinitely many degrees of freedom, but remains in the classical 
framework. In the second case one replaces classical variables by 
quantum variables (operators), but remains within finitely many degree of 
freedoms. Quantum Field Theory inherits the richness of the two 
theories by treating infinitely many quantum variables and enhances them  
further, in particular by interaction, 
particle creation/annihilation, special relativity.

There are thus two paths from the finite-dimensional classical 
calculus to QFT according to the following diagram:
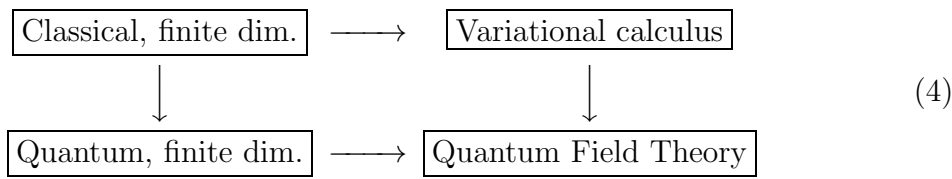
\begin{equation}\label{diag}
\CD
\fbox{Classical, finite dim.}   @>>>\fbox{Variational calculus} 
\\ @VV V	@VV V  \\  \fbox{Quantum, finite dim.}   
@>>> \fbox{Quantum Field Theory} \endCD
\end{equation}
Note now that the passage from ordinary manifolds to variational calculus
did not require a new calculus; for example the notion of derivative 
still make sense replacing points by functions.

On the other hand, the passage from classical to quantum mechanics 
does require a new structure (non-commutativity) a new calculus. 
The standard quantization procedure 
replaces functions by selfadjoint operators and Poisson brackets by 
commutators. In this correspondence $x_h\to P_h$ and
$-i\frac{\partial}{\partial x_h}\to Q_h$ give position and momentum
operators that satisfy the Heisenberg commutation relations 
$[P_h,Q_k]=i\d_{hk}I$.

A quantized, finite-dimensional, calculus has been developed in recent times by A. Connes
(see \cite{C2}); a sample dictionary is here below:
\[
\begin{tabular}{l|l}
CLASSICAL & QUANTUM \\
Variable & Operator\\
Differential & $[F,\cdot]$\\
Integral & $\int\!\!\!\!\! -$ (Dixmier trace)\\
Infinitesimal & Compact operator \\
$\cdots$ & $\cdots$
\end{tabular}
\]
Concerning Quantum Field Theory, more or less implicit suggestions
concerning a ``second quantized'' or QFT calculus can be found in 
\cite{C,D,KK,K,JLO,L2}. In particular we 
consider Jones subfactors and endomorphisms or Connes correspondences to 
be basic objects \cite{L3} in this setting. The underlying structure 
at each level is illustrated in the following table:
\[
\begin{tabular}{|c|l|l|}
\hline
\parbox[l]{2.2cm}{CLASSICAL}&
\parbox[l]{4.5cm}{\centerline{Classical variables}
\par \centerline{Differential forms} \par \centerline{Chern classes}}
   & \parbox[c]{6cm}
   {\centerline{$\phantom{\overset{*}{I}}$Variational 
calculus$\phantom{\overset{*}{I}}$}\par 
 \centerline{Infinite dimensional manifolds}
   \par \centerline{Functions spaces}
   \par \centerline{$\phantom{\underset{*}{I}}$Wiener 
measure$\phantom{\underset{*}{I}}$}}   \\
\hline
\parbox[l]{2.2cm}{QUANTUM}&
\parbox[l]{4.5cm}
{\centerline{$\phantom{\overset{*}{I}}$Quantum 
geometry$\phantom{\overset{*}{I}}$}
\par \centerline{Fredholm operators}
\par \centerline{Index}
\par \centerline{$\phantom{\underset{*}{I}}$Cyclic 
cohomology$\phantom{\underset{*}{I}}$}}
& \parbox[c]{6cm}
{\centerline{Subfactors}
\par \centerline{Correspondences, Endomorphisms}
\par \centerline{Multiplicative index}
\par \centerline{Supersymmetric QFT, $(\gA,\cH,Q)$} }     \\
\hline
\end{tabular} 
\]
Note that there is a non-trivial map from 
\[
\boxed{\text{{\it points}}}\longrightarrow\boxed{\text{{\it fields}}}\ ,
\]
horizontally in the diagram (\ref{diag}), that further enriches the structure. 
At the quantum level this is the second quantization 
functor; this is partly at the basis of the multiplicative structure of 
the index (cf. \cite{B} for an example).

We mention a first result, due to Connes \cite{C}, 
that may be read within the context of QFT analysis: the index map
\[
\text{Ind}_Q : K_0(A(\cO))\to\mathbb Z
\]
is not polynomial and the K-theory group $K_0(A(\cO))$ is of infinite rank. 
Here $A(\cO)$ is a ``smooth'' local Bosonic algebra associated with 
a free massive supersymmetric field on the cylinder.

After these premises, let us discuss
the basic objects of our analysis, superselection sectors.

{\it Superselection sectors as QFT analogs of  elliptic operator}.
The celebrated Atiyah--Singer index theorem equates
the analytic index  of an elliptic operator to a
geometric--topological index. The analytic index is
the Fredholm index, which is manifestly an integer. 
The geometric index  is intrinsically  invariant under deformations.
A major consequence of the index theorem 
is then the integrality of the geometical index.

As discussed, Operator Algebras provide the proper quantization
(non commutative setting) for measure theory, topology and geometry.
In particular, extensions of the index theorem by means of 
noncommutative K-theory
and cyclic cohomology occur naturally in 
Noncommutative Geometry \cite{C}. These results pertain to Connes 
quantized calculus.

Here we shall deal with a quantum field theory
analog of the index theorem. The role of the 
Fredholm linear operators
is now played by the endomorphisms of an infinite 
factor with finite Jones index \cite{L2}. In Quantum Field Theory
the index-statistics theorem \cite{L1} equates the DHR statistical 
dimension with the square root of the Jones index or, 
more precisely, with the minimal dimension, whence
the integrality
of the index is immediate by the integrality of
the statistical dimension \cite{DHR}.  We then look for 
possible  geometric counterparts of the statistical dimension.

Our framework is Quantum Field Theory and its 
superselection structure \cite{WWW},
in the Doplicher-Haag-Roberts framework \cite{DHR}. 
The local observable algebra $\gA(\cO)$
associated with region $\cO$ of the spacetime provides a 
noncommutative version of the algebra of functions with support 
in $\cO$, and the localized endomorphisms with finite statistics 
are analog to the elliptic differential operators, as suggested by S. 
Doplicher \cite{D}. 
Note indeed that an endomorphism $\r$ localized in the double cone $\cO_0$ is local 
in the sense that 
\[
\r(\gA(\cO))\subset \gA(\cO), \quad \cO\supset\cO_0,
\]
similarly to the  locality property that characterizes the
differential operators in the classical 
setting \cite{P}. More is true, the correspondence 
$\cO\to\gA(\cO)$
is endowed with a natural, but not manifest, 
sheaf structure with respect to which covariant localized
endomorphisms are sheaf maps, a fact that will be used
only implicitly (it gives the automatic covariance used 
in Sect. \ref{norm}).

It is now clear that, in our context, geometric information may 
be contained in the 
classical geometry of the spacetime and in the net 
\[
\cO\to\gA(\cO) \ .
\]
A geometrical description of the superselection structure of $\gA$ as 
been given by Roberts \cite{R}, who defined a non-abelian cohomology 
ring $H^1_R(\gA)$ whose elements correspond to the superselection 
sectors of $\gA$. In his formalism, however, it is unclear how to 
integrate a cohomological class to obtain an invariant, the dimension.

A conceptually different cohomological structure appears in \cite{CT} 
by considering unitary cocycles associated with a dynamics. As we 
shall see in Sect. \ref{coho}, if we consider only localized unitary 
cocycles associated with translations, we obtain a cohomology ring 
$H^1_\t(\gA)$ which describes the covariant superselection sectors. 
Denoting by ${\mathfrak S}_{KMS}$ the set of extremal KMS states for 
the time evolution, at inverse temperature $\b$, satisfying Haag duality, 
we have indeed a pairing
\begin{equation}\label{p}
{\mathfrak S}_{KMS}\times H^1_\t(\gA)\ni \f\times [u]\to
  \langle \f,[u]\rangle
=\int\! u(i\b)\text{d}\f \in \mathbb R
\end{equation}
that we shall below describe in equations (\ref{CP},\ref{ac}). We have 
used, only here, the notation $\int\! u\text{d}\f\equiv\f(u)$ in order 
to provide resemblance with the classical context.

It is now useful to compare in a table the context in the Atyiah-Singer 
theorem and in  Quantum Field Theory, including the role possibly 
played by supersymmetry (see below).
\[
\text{{\footnotesize\begin{tabular}{|l|l|l|}
\hline
{}&{\it Atiyah-Singer context} & {\it QFT context}\\
\hline
sheaf structure    & functions on manifold $\cV$   & net (sheaf) of 
C$^*$-algebras on $\cV$\\
smooth structure   & smooth functions    & net of dense $^*$-algebras\\
differential operator & sheaf map              & localized endomorphism\\
elliptic operator  & Fredholm opearator   & finite index endomorphism \\
analitical index   & Fredholm index            & minimal dimension\\
intergrality        & Fredholm index $\in \mathbb Z$ &  statistical 
dimension $\in \mathbb N$\\                                             
geometric index  & associated with $(D,\cV)$ & associated with 
$(\r,\gA,\cV)$ \\
cohomology & De Rham & Roberts; cyclic\\ 
deformation invariance & intrinsic             & perturbation 
invariance\\
Chern character     & Chern character           & pairing (\ref{p}); 
JLO cyclic cocycle\\                                       
Hamiltonian        & Laplacian                 & Killing Hamiltonian\\
grading            & Dirac operator            & supersymmetry\\
spectral formula  & heat kernel               & (super)-Gibbs state\\
\hline
\end{tabular}}}  
\]
In this table we have considered finite volume spaces, although we 
will deal with non-compact spacetimes. As is known, 
Gibbs states becomes KMS states at infinite volume \cite{HHW}. The analogue 
replacement of super-Gibbs functionals by super-KMS functionals 
is far less obvious, see the outlook. However, as mentioned, 
for many purposes concerning the index, one may work with ordinary KMS states. 

{\it Black holes, conformal symmetries and holography.}
There is a first important setting where the above programme may
be implemented, with a formula for the geometric index involving
(classical) spacetime geometry, namely the analysis of charge addition
for a quantum black hole in a thermal state.
 
We introduce here another piece of structure, 
namely we consider Quantum Field Theory on a curved 
spacetime. This physical theory combines General 
Relativity and Quantum Field Theory, but treats the gravitational field
as a background field and therefore disregards effects occurring at the 
Planck length. Yet important effects, as the Hawking effect \cite{Ha},
pertain to this context.

We start with a black hole described by a
globally hyperbolic spacetime with
bifurcate Killing horizon, for example the Schwarzschild-Kruskal
spacetime, and we consider quantum effects on this gravitational 
background. We then consider the incremental entropy due to the
addition of a localized charge. 
The case of the Rindler spacetime has been
previously described by means of a local analogue
of a Kac-Wakimoto formula \cite{L4}.

Here however we use a different point of view and conceptual 
scheme. The first basic point is that the restriction of the net $\gA$
to each of the two horizon components $\mathfrak h_+$ and  
$\mathfrak h_-$
gives a conformal net on $S^1$, a general fact that is obtained
by applying Wiesbrock's characterization of conformal nets \cite{W},
a structure already discussed in \cite{L5,GLRV}. 
It may appear analogous to the holography on the anti-de 
Sitter spacetime that independently appeared in 
the Maldacena-Witten conjecture\footnote{We thank M. 
Konsevich for pointing out this similarity to 
us.} \cite{Ma,Wi}, proved by Rehren \cite{Re}. Yet their context 
differs inasmuch as
the anti-de Sitter spacetime is not globally hyperbolic and 
the holography there is a peculiarity of that spacetime, rather than
a general phenomenal.

With these conformal nets at one hand, and assuming the due duality 
properties,
we may consider endomorphisms $\r$ and $\s$ of $\gA$ that are
localizable on $\mathfrak h_+$ or  $\mathfrak h_-$. 
In particular, we will have the formula
for the difference of the logarithm of dimensions:
\[
\log d(\r)-\log d(\s)=
{\frac{2\pi}{\k(\cV)}}(F(\f_\r|\f_\s)+F(\f_{\bar\r}|\f_{\bar\s})),
\]
where $F$-terms represent the incremental free energy between
the thermal equilibrium states with the charges $\r$ and $\s$ or the conjugate
charges $\bar\r$ and $\bar\s$. Here $\f$ is the Hartle-Hawking state 
and $\f_\r$ is the corresponding equilibrium state in presence of the 
charge $\r$. The geometry appears here in the surface
gravity $\k(\cV)$ associated to the spacetime manifold $\cV$,
see \cite{KW}.
As a consequence the right hand side is the difference of the
logarithm of two integers.

We will consider also different temperature states. 
For a finer analysis of this context we refer to our 
Sect. \ref{blackhole}. 
We shall need to study the chemical potential, 
as we are going to explain. 

{\it Chemical potential. Relativistic case.}
The chemical potential  is a label that each charge sets 
on different equilibrium states at 
the same temperature. Its structure in Quantum 
Statistical Mechanics has been explained in the work of
Araki, Haag, Kastler and Takesaki \cite{AHKT}, see also \cite{AK}. 
The labels appear by considering the extensions of these states 
from the observable algebra to KMS states of the field algebra 
(with the time evolution modified by one-parameter subgroups 
of the gauge group).

Our aim is to consider temperature states in Quantum Field Theory.
This may be motivated by wish to study extreme physical contexts such
as the early universe (not discussed in this paper) or black holes. 

In the context of Quantum Field Theory on the $d+1$-dimensional
Minkowski spacetime with $d\geq 2$, Doplicher and Roberts \cite{DR}
have constructed the field algebra associated with local observables
and short range interaction charges.

It is not difficult to apply the AHKT analysis,
originally made in the context of C$^*$-dynamics, to the case
of Quantum Field Theory on the Minkowski spacetime as above.
It turns out the every covariant irreducible localized endomorphism $\r$ with
finite dimension extends to the weak closure in the GNS 
representation $\pi_\f$ associated with a KMS state $\f$, a fact that
should be expected on physical grounds because the addition of a 
single charge should not lead to an inequivalent representation
for an infinite system.

If $\f$ is extremal KMS, we are then led to the context of 
endomorphisms of factors with finite index. Assuming 
$\pi_\f$ to satisfy Haag duality, the extension of $\r$ to the weak 
closure is still irreducible.
If $u$ is the time 
covariance cocycle for $\r$, 
then $u$ is a Connes Radon-Nikodym cocycle \cite{C1} up to phase, hence
it satisfies certain holomorphic properties. 
We have the formula, that generalizes \cite{L4},
\begin{equation}\label{CP}
\log d_\f(u) = \log d(\r) + \b\mu_\r(\f)\ ,
\end{equation}
where the holomorphic dimension is defined by 
\begin{equation}\label{ac}
d_\f(u)\equiv\acb \f(u(t))\ .
\end{equation}
If there is a canonical choice for $u$, for example if the $\r$ is 
Poincar\'e covariant in the vacuum sector,
then the above formulae define the chemical potential 
$\mu_\r(\f)$ of $\f$ corresponding to 
the charge $\r$. This extends an analogous expression in
the AHKT work in the case of abelian charges ($d(\r)=1$).
In general only the difference $\mu_\r(\f) - \mu_\r(\psi)$
is intrinsic.

Now the localized endomorphisms form a C$^*$-tensor category 
\cite{DHR} and in a C$^*$-tensor category there exists a natural 
anti-linear conjugation
on the arrows between finite-dimensional objects \cite{LR}:
\[
T\in(\r,\s)\to T^{\bullet}\in(\bar\r,\bar\s)\ .
\]
Therefore, even if in general $u$ is defined only up to phase,
$u^{\bullet}$ is a covariance cocycle for the conjugate charge 
$\bar\r$ with opposite chemical potential 
\begin{equation}\label{asym}
\mu_{\bar\r}(\f)=-\mu_\r(\f)
\end{equation}
and we obtain an expression
for the intrinsic dimension $d(\r)$ as the geometric mean
\[
d(\r)= \sqrt{d_\f(u)d_\f(u^{\bullet})}\ ,
\]
which is independent of the phase fixing.
We regard the right hand side of this expression as a geometric 
dimension, according to what was explained before, being candidate for 
geometric interpretation.

The quantity $\b^{-1}\log d_\f(u)$ represents the incremental free 
energy (adding the charge $\r$). It has then a canonical decomposition
as a sum of an intrinsic part $\b^{-1}\log d(\r)$, which is
independent of $\f$, where $\log d(\r)$ is half the 
incremental entropy associated with $\r$ (cf. \cite{PP,L4}), 
and the chemical potential part $\mu_\r(\f)$ 
characterized by the asymmetry with respect to charge conjugation 
in (\ref{asym}).

{\it Chemical potential. Low dimensional case.}
Motivated by black hole thermodynamics and the associated
conformal nets on the black hole horizon,
among other considerations, one is led to the analysis of
the chemical potential in one-dimensional Quantum Field Theory.
In this context the field algebra does not any longer exist and the
AHKT work is not applicable.

The notion of holomorphic dimension still makes sense and is the
basis of our analysis. The crucial point however is to show that 
localized endomorphisms are normal in the representation $\pi_\f$
associated to a thermal state $\f$. Based on Wiesbrock 
characterization of conformal nets on $S^1$ \cite{W}, we shall see 
that there is a conformal net associated with $\f$, the conformal
thermal completion. If $\pi_\f$ satisfies duality 
for half-lines, we prove that the thermal completion 
automatically satisfies Haag 
duality, namely it is strongly additive \cite{GLW}. 

Then, a localized endomorphism $\r$ 
with finite dimension of the original net
gives rise to a transportable localized endomorphism of the thermal 
completion. By a result in \cite{GL1} $\r$ is automatically 
conformally covariant. We finally use this conformal covariance
to show the normality of $\r$ in the representation $\pi_\f$.

For completeness, we extend the work of AHKT also in regard 
to describing the chemical potential in terms of extensions 
of KMS states. This will be achieved by considering extensions
to the quantum double, a  C$^*$-algebraic version of a 
construction in \cite{LR}, that is a substitute for the non-existing
field algebra.

At this point the analysis of the chemical potential goes through as
in the higher dimensional case with the corresponding formulae
for the incremental free energy. In the case of globally hyperbolic
spacetimes with bifurcate Killing horizon, these results allow one
to treat thermal states for the Killing evolution and charges
localizable on the horizon, as explained.

{\it The expected role of supersymmetry.}
Connes cyclic cohomology enters in
Supersymmetric Quantum Field Theory via the work of
Jaffe, Lesnieski and Osterwalder \cite{JLO}, see also \cite{K}.
There is a noncommutative 
Chern character associated with a thermal equilibrium state, 
i.e. an entire cyclic 
cocycle associated to a supersymmetric KMS functional for 
the time evolution. In this context an index formula appears, 
which essentially coincides with our
previous formula (\ref{ac}).  There, however, the index formula 
acquires the geometric meaning of evaluating JLO cycle at 
the identity. Indeed such a formula was used
to show the deformation invariance of the index.

On the other hand, in our context,
our formula for the dimension is independent of an underlying
supersymmetric structure. When we consider unitary
cocycles associated with charges localizable in a bounded region,
the dimension varies, but remains an integer. 

If we consider localized endomorphisms and super-KMS functional,
namely if we consider the index associated with a covariance
cocycle and a super-KMS functional,
we should expect our index formula to be read geometrically
by the JLO cyclic cocycle. We will explain this point at the end of 
this paper, as it can certainly give further insight. But the present 
picture is too primitive to be directly applicable because
super-KMS functionals fail to exist  when the spacetime is 
non-compact in the most natural situation when
the functional is translation invariant and space translations act
in an asymptotically abelian fashion \cite{BL}.
This drawback is entirely caused by the assumption that the super-KMS
functional is bounded. The structure associated with unbounded
super-KMS functional is presently under investigation.
                                                           
\section{First properties of holomorphic cocycles.}
In this section we begin to study holomorphic cocycles and give
first formulae for the dimension.
\subsection{Index formulae. Factor case.}\label{factorcase}
Here we recall and extend results in \cite{L4}, Section 1. We
shall show how to obtain a holomorphic formula for the dimension 
of a sector that  assumes neither a perfect symmetry group nor a PCT 
anti-automorphism. The reader is however assumed to have read
the above quoted reference.

Let $\cM$ be an infinite factor and denote by End$(\cM)$ the 
(injective, normal, unital)\footnote{The `index' or 
`dimension' terminology is here interchangeable; 
analogously, the Fredholm index of an
isometry is the dimension of its cokernel.} 
endomorphisms of $\cM$ with 
finite Jones index and Sect$(\cM)$ the sectors 
of $\cM$, namely the equivalence classes of End$(\cM)$ modulo 
inner automorphisms of $\cM$.
End$(\cM)$ is a tensor
category where the tensor product is the composition, cf. \cite{DR}.
The intertwiner space $(\r,\s)$ between objects $\r, 
\s\in{\rm End}(\cM)$ is defined as
\[
(\r,\s)\equiv\{ T\in\cM: T\r(X)=\s(X)T, \forall\ X\in\cM\}\ .
\]
Sect$(\cM)$ is endowed with a natural conjugation. 
$\bar\r\in{\rm End}(\cM)$ is a conjugate of $\rho$ iff the conjugate 
equation holds true: there exist multiples of isometries 
$R\in(\iota,\bar\r \r)$ and $\bar R\in(\iota,\r\bar\r)$, 
that we normalize with $||R||=||\bar R||$, such that
\begin{equation}
	R^*\bar\r(\bar R)=1,\quad
	\bar R^*\r(R)=1\  .
	\label{}
\end{equation}
The minimum
\[
d(\r)\equiv {\rm min}\|R\| \|\bar R\|
\]
over all possible choices of $R$ and $\bar R$ is the
{\it intrinsic dimension} of $\r$. A pair $R_\r$, $\bar R_\r$ 
where the minimum is attained exists and coincides with a standard solution
as discussed in \cite{LRo}.
It turns out that 
$d(\r)=d_{an}(\r)$, the analytical dimension defined as 
$d_{an}(\r)=\sqrt{[\cM:\r(\cM)]}$, the square root
of the minimal index  $[\cM:\r(\cM)]$  
(Jones-Kosaki index with respect to the minimal expectation) \cite{L2}.

For each $T\in (\r_1,\r_2)$, the conjugate arrow 
$T^{\bullet}\in (\bar\r_1,\bar\r_2)$ is  defined by  
\[\label{bullet}
T^{\bullet}=\bar\r_2(\bar R_{\r_1}^*T^*)R_{\r_2}
\]
where $R_{\r_i}$ and $\bar R_{\r_i}$ give a standard solution for
the conjugate equation defining the conjugate $\bar\r_i$.
The map $T\to T^{\bullet}$ is anti-linear and satisfies natural 
properties, see \cite{LRo}.

Fix a normal faithful state $\f$ of $\cM$ and let $\s^\f$ be the modular 
group of $\f$. As is well known, $\s^\f$ satisfies the KMS condition 
with respect to $\f$ at inverse temperature $-1$, namely, setting 
$\a_t=\s^\f_{-t}$, the relation (\ref{KMS}) holds 
with ${\b=1}$. For a fixed $\r$, let $u(\r,\cdot)$ be a unitary 
$\s^\f$-cocycle (i.e. 
$u(\r,t)$ is a unitary in $\cM$ and 
$u(\r,t+s)=u(\r,t)\s^\f_t(u(\r,s))$, $t,s\in\Reali$ ) such that
\[
{\rm Ad}u(\r,t)\cdot\s_t^\f\cdot\r =\r\cdot\s_t^\f\ ,
\]
that is $u(\r,t)\in(\r_t,\r)$, where 
$\r_t\equiv \s^\f_t\cdot\r\cdot\s^\f_{-t}$.
Note that $u(\r,\cdot)$ is not not assumed to be continuous.
Once standard $R_\r$ and $\bar R_\r$ are given, 
we assume the corresponding 
operators for the conjugate equation for $\r_t$ and $\bar\r_t\equiv
\s^\f_t\cdot\bar\r\cdot\s^\f_{-t}$ to be given by 
$R_{\r_t}=\s_t^\f(R_\r)$ and $\bar R_{\r_t}=\s_t^\f(\bar R_\r)$.

Then $u(\r,t)^{\bullet}\in(\bar\r_t,\bar\r)$ is given by
\[
u(\r,t)^{\bullet}=\bar\r(\bar R^*_{\r_t}u(\r,t)^*)R_\r
= \bar\r(\s^\f_t(\bar R^*_\r) u(\r,t)^*)R_\r \ .
\]
If $\r$ is irreducible, the choice of
$R_\r$ and $\bar R_\r$ is unique up to a phase,
therefore $u(\r,t)^{\bullet}$ is uniquely defined. This holds in more
generality, if $\r$ is reducible.
\begin{Prop}\label{1.2}  $u(\r,t)^{\bullet}$ is well defined, namely it does not
depend on the choice of  $R_\r$ and $\bar R_\r$ giving   
a solution for the conjugate equation for $\r$ and $\bar\r$.
\end{Prop}
\begin{proof}
Let $R_\r$ and $\bar R_\r$ be a standard solution.
If $R'_\r$ and $\bar R'_\r$ is another solution of the 
conjugate equation, then $R'_\r= \bar\r(v)R_\r$ and 
$\bar R'_\r = v^{*-1}\bar R_\r$ for some invertible $v\in(\r,\r)$
\cite{LRo}, hence the conjugate of $u(\r,t)$ with respect to
$R'_\r$ and $\bar R'_\r$ is given by
\begin{multline}
\bar\r({\s^\f}_t (\bar R^*_\r) \s^\f_t(v^{-1})u(\r,t)^*)\bar\r(v)R_\r 
= \bar\r(\s^\f_t (R^*_\r) u(\r,t)^*\s^{\f_\r}_t(v^{-1}))\bar\r(v)R_\r \\
= \bar\r(\s^\f_t (R^*_\r) u(\r,t)^* v^{-1})\bar\r(v)R_\r 
=\bar\r(\s_t ^\f({\bar R_\r}{^*})) u(\r,t)^*)R_\r = u(\r,t)^{\bullet}\ ,
\end{multline}
where $\s^{\f_\r}$ is the modular group of $\f\cdot\Phi_\r$ , so that
$\s^{\f_\r}_t(v) = v$ because the minimal left inverse $\Phi_\r$ 
of $\r$ is tracial on $(\bar\r,\bar\r)$.
\end{proof}
\begin{Prop}\label{geo} Let $\r$ be irreducible. We have 
\[
d(\r)^2=\aci \f(u(\r,t))\f(u(\r,t)^{\bullet})\ .
\] 

If $u(\r,\cdot)$ is weakly continuous, then $u(\r,\cdot)$ and 
$u(\r,\cdot)^{\bullet}$ are 
holomorphic (see 
below) in the state $\f$, and 
\[
d(\r)^2 = \aci\f(u(\r,t))\aci\f(u(\r,t)^{\bullet})\ .
\]
\end{Prop}
\begin{proof} We may suppose that $d(\r)<\infty$. Set $\f_{\r}\equiv
\f\cdot\Phi_{\r}$, where $\Phi_{\r}$ is the minimal left inverse of
$\r$ as before. 
Since the Connes Radon-Nikodym cocycle $(D\f_{\r}:D\f)$ is
a covariance cocycle for $\r$ \cite{L4}, 
there exists a one-dimensional character $\chi$ of $\Reali$ 
such that
\[
u(\r,t)=\chi(t)d(\r)^{it}(D\f_{\r}:D\f)_t\ ,
\]
(both cocycles intertwine $\s^{\f}$ and $\s^{\f_{\r}}$).
As shown in \cite{L4}, $(d(\r)^{it}(D\f_{\r}:D\f)_t)^{\bullet}
=d(\r)^{it}(D\f_{\bar\r}:D\f)_t$
\[
u(\r,t)^{\bullet}=\overline{\chi(t)}(d(\r)^{it}(D\f_{\r}:D\f)_t)^{\bullet}
=\overline{\chi(t)}d(\r)^{it}(D\f_{\bar\r}:D\f)_t\ ,
\]
Moreover the Connes cocycle is holomorphic and 
$\ac \f((D\f_{\r}:D\f)_t) =\f_\r(1)=1$ \cite{C}, therefore
\begin{multline}
\aci\f(u(\r,t))\f(u(\r,t)^{\bullet})\\
=\aci d(\r)^{2it}\f((D\f_{\r}:D\f)_t)
\f((D\f_{\bar\r}:D\f)_t) = d(\r)^2\ .
\end{multline}
If $u(\r,\cdot)$ is continuous, then $\chi$ is continuous. So 
$\chi$ extends to an entire function,  hence 
$u(\r,\cdot)$ is holomorphic, and the second formula
in the statement follows from the first one.
\end{proof}
{\it Remark.} Let $\cT\subset\text{End}(\cM)$ be a C$^*$ tensor category
as before, and $u(\r,t)$ a two-variable cocycle. Setting $d_\f(\r)
=d(u_\r)\equiv \aci\f(u_\r(t))$ (see below) the arguments in \cite{L4} 
show that
\begin{align*}
d_\f(\r_1\oplus\r_2)&=d_\f(\r_1)+d_\f(\r_2),\\
d_\f(\r_1 \r_2)&=d_\f(\r_1)d_\f(\r_2),\quad \r_1 , \r_2\in\cT\ .
\end{align*}
In particular, if $\cT$ is rational, namely there are only finitely 
many inequivalent irreducible objects, the usual application of the 
Perron-Frobenious theorem entails $d(\r)=d_\f(\r)$ for all objects of 
$\cT$ (no chemical potential, see Sect. \ref{QFTM}). 
\subsubsection{Case of a non-full C$^*$ tensor sub-category.}
Now let $\cT$ be a C$^*$ tensor category with conjugates contained
in End$(\cM)$, thus the objects of $\cT$ are finite-index endomorphisms of
$\cM$, but we do not assume that $\cT$ is a full sub-category of 
End$(\cM)$, namely the intertwiner spaces $(\r,\s)$ in $\cT$ can be
strictly contained in the corresponding intertwiner spaces in End$(\cM)$.

Assume that the modular group $\s^\f$ gives an action of $\Reali$ on 
$\cT$, that is $\r_t\in\cT$, for all $t\in \Reali$,  $\r\in\cT$, 
and $\s^{\f}_t((\rho,\rho'))=(\r_t,\r'_t)$ if $\r,\r'\in\cT$ are objects 
of $\cT$. 

Recall that $u$ is a {\it two-variable unitary cocycle} for the above action
if, for each fixed object $\r\in\cT$, $u(\r,\cdot)$ is a unitary 
$\s^\f$-cocycle as above and, for each fixed $t\in\Reali$, $u(\r,\cdot)^*$ is 
cocycle with respect to $\s^\f_t$, namely
\[
u(\r\s,t)=\r(u(\s,t))u(\r,t), \quad \r,\s\in\cT,
\]
and
\[
Tu(\r,t)=u(\s,t)\s^\f_t(T),\quad \r,\s\in\cT,\ T\in (\r,\s).
\]
Note that if $u$ is a two variable cocycle also for  
the full tensor subcategory of End$(\cM)$ with the same objects of 
$\cT$, then $u(\r,t)^{\bullet}$ defined there coincides with
$u(\r,t)^{\bullet}$ defined in $\cT$ cf. \cite{L4}, Propositions
1.5 and A.2.
\begin{Cor}\label{2v}With the above notations, if $u(\r,t)$ is a 
weakly continuous unitary two-variable 
cocycle for the action of $\Reali$ on $\cT$ given by $\s^\f$, then 
\[
d(\r) \leq\sqrt{\aci\f(u(\r,t))\f(u(\bar\r,t))}
\]
for all irreducible $\r\in\cT$\ . Here $d(\r)$ is the intrinsic dimension of 
$\r$ as an object of $\cT$.
\end{Cor}
\begin{proof} We have $u(\r,t)^{\bullet}=u(\bar\r,t)$ \cite{L4}
(the conjugate map is the one associated with $\cT$), 
hence the above Proposition \ref{1.2} applies, provided $\r$ is irreducible
in End$(\cM)$. If  $\r$ is reducible in End$(\cM)$
we have $u(\r,t) = z(t)d_{an}(\r)^{it}(D\f_\r : D\f)_t$ with $z(t)\in(\r,\r)$. By 
using the tracial property of $\Phi_\r$ it is easy to check that
$z$ is a one-parameter group of unitaries in the finite-dimensional 
algebra $(\r,\r)$, and therefore can be diagonalized. If 
$\r=\oplus_i \r_i$ is a decomposition of $\r$ into irreducibles 
(with eigen-projections of $z$) in 
End$(\cM)$, we have a corresponding decomposition of
$d_{an}(\r)^{it}(D\f_\r : D\f)_t$ as direct sum of the 
$d_{an}(\r_i)^{it}(D\f_{\r_i} : D\f)_t$ thus
$d_\f(u(\r_i,\cdot))=\ell_i d_{an}(\r_i)$ with $\ell_i>0$, hence
\[
d_\f(u_\r)d_\f(u_{\bar\r} )=\sum_{i,j}\ell_i d_{an}(\r_i)
\ell_j d_{an}(\r_j)\geq \sum_{i,j}d_{an}(\r_i)
d_{an}(\r_j)=d_{an}(\r)^2\geq d(\r)^2\ .
\]
\end{proof}
We now recall that the following holds.
\begin{Prop}{\rm \cite{L4}.}\label{PCT} In the setting of Prop. \ref{1.2}, 
if there exists an 
$\f$-preserving anti-automorphism $j$ of $\cM$ inducing an 
anti-automorphism of $\cT$ such that $j\cdot\r\cdot j=\bar\r$ and
$j(u(\r,t))=u(\bar\r,-t)$, then $d(\r)=\aci \f(u(\r,t))$ for all 
objects $\r\in\cT$.
\end{Prop}
\begin{proof} See \cite{L4}, Prop.1.7.
\end{proof}
\subsection{The holomorphic dimension in the C$^*$-case.}
\label{sec:first}
In this section we give a first look at the structure that emerges 
in the C$^*$ context, in analogy to what studied in the previous section in 
the setting of factors. Here we assume from the start 
that a holomorphic dimension is definable, 
postponing the more relevant derivation of the holomorphic 
property and the analysis of the chemical potential to subsequent 
sections.

Let $\gA$ be a unital C$^*$-algebra and $\a$ a
one-parameter automorphism group  of $\gA$.
A linear functional $\f\in \gA^*$ is said 
to be a {\it  KMS functional}
with respect to $\a$ at inverse temperature $\b>0$ 
if for any given $a,b\in \gA$ there exists a 
function $F_{a,b}\in A(S_\b)$ such that
\begin{align}\label{KMS}
	F_{a,b}(t) & = \f(\a_t(a)b)\\
	F_{a,b}(t+i\b) & =  \f(b\a_t(a))
\end{align}
Here $S_\b$ is the strip $\{0<\mbox{Im}z < \b\}$ and $A(S_\b)$
is the algebra of functions
analytic in $S_\b$, bounded and continuous on the closure of $S_\b$.
We do not assume $\a$  to be pointwise norm  continuous, nonetheless
a weaker continuity property follows from the KMS condition.
Note that a KMS functional $\f$ is $\a$-invariant.

Let now $u\in\gA$ be a unitary cocycle with respect to $\a$, 
namely $t\in\Reali\to u(t)\in \gA$ is a map taking values in the 
unitaries of $\gA$ satisfying the equation
$$
u(t+s)=u(t)\a_t(u(s)).
$$
We shall say that the cocycle $u$ is {\it holomorphic}, in the 
functional $\f$, if the function 
$t\in\Reali\to \f(u(t))$ is the boundary value of a function in 
$A(S_\b)$.

If $u$ is holomorphic in the  state $\f$,
we define the {\it holomorphic dimension}
of the cocycle $u$ (with respect to $\f$) by
\[
\label{holdim}
d_\f(u) =\acb\f(u(t)).
\]
As we shall see, in our context 
$d_\f(u)$ will be a positive number related to
a  (noncommutative) relative index.

Clearly, for a given dynamics $\a$, $d_\f(u)$  
may depend on the KMS state
$\f$. We shall sometimes rescale the ``time parameter'' to make
the inverse temperature $\b=1$. If $u$ is not holomorphic we write 
$d_\f(u)=+\infty$.
 
Let $\r$ be  an endomorphism of $\gA$.
We shall say that $\r$ is {\it covariant} if there exists a  
$\a$-cocycle of unitaries $u(\r,t)\in \gA$  such that 
\begin{equation}\label{cov}
\a_t\cdot\r\cdot\a_{-t} =\text{Ad}u(\r,t)^*\cdot\r.
\end{equation}
We shall say that $\r$ has {\it finite holomorphic dimension} (with respect to
the KMS state $\f$) if it is 
covariant and there exists a covariance cocycle 
$u(\r,\cdot)$ as above with finite holomorphic dimension. 
Note that, if $\gA$ has trivial 
centre and $\r$ is irreducible, 
$u(\r,\cdot)$ is unique up to a phase, that doesn't alter the
finite-dimensional property of $u(\r,\cdot)$, provided such a phase 
is chosen to be a continuous character of $\Reali$.

In the rest of this section we study whether endomorphisms of $\gA$ with 
finite holomorphic dimension extend to $\pi_\f(\gA)''$.

In the following we identify $\gA$ with 
its image $\pi_\f(\gA)$ and suppress the suffix $\f$.

\begin{Lemma}\label{ext}
Let $\gA$ be a C$^*$-algebra acting on a Hilbert space $\cH$ and $\r$ an 
endomorphism of $\gA$. Let $\x$ be a cyclic separating vector for 
$\cM=\gA''$ and $\f=(\cdot\, \x,\x)$. If $\f\cdot\r$ extends to a normal 
faithful positive functional of $\cM$,
then $\r$ extends to a normal endomorphism of $\cM$.
\end{Lemma}
\begin{proof}
Let $\eta$ be a cyclic separating vector such that $(\r(a)\x,\x)
= (a\eta,\eta)$,
 $a\in\gA$,
and let $V$ be the isometry of $\cH$ given by
 \begin{equation}
 	Va\eta=\r(a)\x, \quad a\in \gA\ .
 	\label{}
 \end{equation} 
The final projection  of $V$ is given by
\begin{equation}
	e=VV^{*}=\overline{\r(\gA)\x}\in\r(\gA)'
	\label{}
\end{equation}
thus $x\in \cM\to VxV^*$ is a homomorphism of $\cM$ onto $\r(\gA)''e$
and 
\begin{equation}
	VaV^*=\r(a)e, \quad a\in \gA.
	\label{}
\end{equation}
Now the central support of $e$ in $\r(\gA)''$ is $1$ as 
$\overline{\r(\gA)'\x}\supset \overline{\cM'\x}=\cH$, hence if $x\in \cM$ 
there exists a unique $\r(x)\in\r(\gA)''$ such that $\r(x)e = VxV^*$,
providing an extension of $\r$ to $\cM$. As $\eta$ is separating, $\r$ is 
an isomorphism.
\end{proof}
Now, as in the factor case, End$(\gA)$ is the tensor category 
whose objects $\r,\s,\dots$ are the
endomorphisms of $\gA$: the monoidal product $\r\otimes\s =\r\s$ is given
by the composition of maps, while the intertwiner space $(\r,\s)$ is 
given by $\{T\in\gA : T\r(a)=\s(a)T,\, \forall\; a\in\gA\}$. The tensor 
product of intertwiners is also defined in a natural fashion, see e.g.
\cite{LRo}. 
The conjugate equation is defined as in the previous section. The 
\emph{intrinsic dimension} $d(\r)$, and the conjugation on arrows 
are defined as well, in fact all these notions make
sense for a a general tensor C$^*$-category,  \cite{LRo}.
The following proposition is a special case of results in \cite{LRo}. We 
state it in the particular case needed for our applications.
\begin{Prop}\label{extbraided}{\rm (\cite{LRo}).}
Let $\gA$ be a unital C$^*$-algebra, acting on a Hilbert space,
with $\cM=\gA''$ a factor. Let $\cT$ be a tensor category with 
conjugates and subobjects of endomorphisms of $\gA$ admitting
a unitary braid group symmetry. Suppose that every endomorphism
$\r\in\cT$ extends to a normal endomorphism of $\hat\r\in\cM$.
Then 
\[
d(\hat\r)=d(\r)
\]
where $d(\hat\r)=d_{an}(\hat\r)$, i.e. $d(\hat\r)^2$ is the minimal 
index $[\cM:\hat\r(\cM)]$, and $d(\r)$ is the intrinsic dimension 
of $\r$ in $\cT$.
\end{Prop}
\begin{proof}
The map $\r\to\hat\r$ is a functor of C$^*$ tensor categories from 
$\cT$ to a sub-tensor category of End$(\cM)$, hence $d(\hat\r)\leq
d(\r)$. As $\cT$ has a unitary braiding, every real object $\s\in\cT$
is amenable \cite{LRo} Th. 5.31, thus $d(\s)=||m^\s||$, where $||m^\s||$ is the $\ell^2$ 
norm of the fusion matrix $m^\s$ associated with $\s$, therefore
$d(\s)=||m^\s||\leq ||m^{\hat\s}||\leq d(\hat\s)\leq d(\s)$, thus 
$d(\hat\s)=d(\s)$. For any $\r\in\cT$, the object $\s=\r\bar\r$ is
real hence, by the multiplicativity of the dimension, $d(\hat\r)=d(\r)$. 
\end{proof}
\subsubsection{Case of a unique KMS state.}
We now restrict our attention to the case of a unique KMS functional. 
This is done more with an illustrative intent, 
rather than for later applications, where we shall treat a more general
context.
\begin{Prop}\label{extension} Suppose $\f$ is the unique KMS 
functional for $\a$. Then a covariant endomorphism $\r$ of $\gA$ with 
finite holomorphic dimension has a normal extension 
to $\cM\equiv\gA''$.
\end{Prop}
\begin{proof} Let $u=u(\r,\cdot)$ a holomorphic unitary covariance 
$\a$-cocycle.
As $t\to \f(u(t))$ is continuous, the map $t\in\Reali\to u(t)\in \cM$ is 
strongly continuous. To check this, note that by cocycle property it is
enough to verify the continuity at $t=0$ because then the strong limit
\begin{equation}
u(s+t)=u(s)\a_s(u(t))\to 0, \quad \text{as}\; t\to 0
\end{equation}
due to the normality of $\a_s$. 

Let then $x\in \cM$ be a weak limit point of $u(t)$ as $t\to 0$. Then 
$\|x\|\leq 1$ and
\begin{equation}
(x\x_\f,\x_\f)=\lim_{i\to\infty}(u(t_i)\x_\f,\x_\f)=
\lim_{i\to\infty}\f(u(t_i))=\f(1)=(\x_\f,\x_\f)
\end{equation}
for some sequence $t_i\to 0$. Thus $x\x_\f=\x_\f$ by the limit case 
of the Schwartz inequality, thus $x=1$ because $\x_\f$ is separating.

Therefore by Connes' theorem \cite{C} there exists a normal faithful semifinite
weight $\f_{\r}$ on $\cM$ with $(D\f_{\r},D\f)_t = u(t)$ and, by the 
finite holomorphic dimension assumption, $\f_\r(1)=d_\f(u)\leq\infty$ 
so that $\f$ is indeed a positive linear functional.
Then $\f_\r$ is a KMS functional for its modular group 
$\a^{\r}_{-t}=\text{Ad}u(-t)\cdot\a_{-t}$ (we are setting $\b= 1$ here).

The functional on $\gA$ 
\begin{equation}\label{o'}
	\f'(a)=\f_\r(\r(a)),\, a\in \gA
\end{equation}
is KMS with respect to $\a$
\begin{multline}
\ac \f'(a\a_t(b))=\ac \f_\r(\r(a)\r(\a_t(b)))\\
=\ac \f_\r(\r(a)\a^{\r}_t(\r(b)))=\f_\r(\r(b)\r(a))=\f'(ba).
\end{multline}

Hence, by the uniqueness of the KMS state, 
\begin{equation}
	\f'=\l\f
\end{equation}
on $\gA$, for some $\l>0$. Thus $\f_\r$ is a normal faithful functional
and $\f_\r\cdot\r$ is normal too. Therefore
the proof is completed by Lemma \ref{ext}.
\end{proof}
To shorten notation, we shall often set
\[
u_\r =u(\r,\cdot)\ .
\]
\begin{Prop}\label{extensionconj} Let $\f$ be the unique KMS state 
as in Proposition \ref{extension}.
If $\r$ and $\bar\r$ are conjugate and $d_\f(u_\r)<\infty$,
$d_\f(u_{\bar\r})<\infty$, then the extension of $\r$ 
to $\cM$ has finite Jones index.
\end{Prop}
\begin{proof}
By assumption and Proposition \ref{extension} both $\rho$ and 
$\bar\rho$ extend to $\cM$. As $R_\r$ and $\bar R_\r$ are also 
intertwiners on $\cM$ by weak continuity and the conjugate equation 
for $\r$ and $\bar\r$ is obviously satisfied on $\cM$, the extension
of $\r$ to $\cM$ has finite index.
\end{proof}
\begin{Prop}\label{dim}
Suppose that $\f$ is faithful and the unique KMS state for 
$\a$, as in Prop. \ref{extensionconj}. Let $\cT$ be a tensor category 
with conjugates of endomorphisms of $\gA$
and $u(\r,t)\in\gA$ a unitary two-variable cocycle.
If $u_\r $ has finite holomorphic dimension for all irreducible 
objects $\r$, then the intrinsic dimension $d(\r)$ of $\r$ 
is bounded by
\begin{equation}
	d(\r)\leq\sqrt{d_\f(u_\r )d_\f(u_{\bar\r} )}.
	\label{}
\end{equation}
and equality holds if $\r$ is irreducible and extending to $\cM$ 
is a full functor (thus $\r$ is irreducible on $\cM$).
\end{Prop}
\begin{proof} By Lemma \ref{extensionconj} $\r$ extends to $\cM$ and has
finite index. We are then in the case covered by Prop. \ref{geo} 
and Corollary \ref{2v}.
\end{proof}
Motivated by the above Proposition, we define the {\it 
geometric dimension}
$d_{geo}(\r)$ as
\[
d_{geo}(\r)\equiv\sqrt{d_\f(u_\r )d_\f(u_\r ^{\bullet})}\ .
\]
If $\r$ is irreducible, $d_{geo}(\r)$ does not depend on
the choice of the covariance cocycle $u(\r,t)$, because, if we
multiply $u(\r,t)$ by a phase $\chi(t)$, then $u(\r,t)^{\bullet}$
has to be replaced by $\overline{\chi(t)}u(\r,t)^{\bullet}$. Of course,
since for a two-variable cocycle $u(\r,t)$ we have 
\[
u(\r,t)^{\bullet}=u(\bar\r,t)\ ,
\]
in this case we also have
\[
d_{geo}(\r)=\sqrt{d_\f(u_\r )d_\f(u_{\bar\r} )}\ .
\]
Also, since $u(\r\bar\r,t)=\r(u(\bar\r,t))u(\r,t)$, we have
\[
d_{geo}(\r)=\sqrt{d_\f(u_{\r\bar\r})}
\]
if $\r$ and $\bar\r$ extend to $\cM$.

A priori $d_{geo}(\r)$ might depend on the KMS functional $\f$,
but, as we shall see, in most interesting cases it will actually be 
independent of $\f$.
\subsubsection{Graded KMS functionals: 
reduction to ordinary KMS states.}
The above results extend to graded KMS functionals. Indeed the 
analysis of these functionals can be reduced to the case of 
ordinary KMS states.

Let $\gA$ be a $\Interi_2$-graded unital C$^*$-algebra,
namely $\gA$ is a unital C$^*$-algebra equipped with  an involutive 
automorphism $\g$ \footnote{In this paper morphisms always commute 
with the $^*$-mapping and preserve the unit.}. 

Given a graded one-parameter automorphism group $\a$ of $\gA$ 
(i.e. one commuting with $\g$), a linear functional $\f\in \gA^*$ 
is said to be a {\it graded KMS functional}
with respect to $\a$ at inverse temperature $\b>0$ 
if for any given $a,b\in \gA$ there exists a 
function $F_{a,b}\in A(S_\b)$ such that
\begin{align}\label{GKMS}
	F_{a,b}(t) & = \f(\a_t(a)b)\\
	F_{a,b}(t+i\b) & =  \f(\g(b) \a_t(a))
\end{align}
Note that a graded KMS functional $\f$ is $\g$-invariant.

We recall the following.
\begin{Prop}\label{BL}{\rm (\cite{St,BL})}
Let $\f$ be a graded KMS functional of $\gA$ for $\a$ and 
let $\o=|\f|$ be the modulus of $\f$. Then
$\o$ is an ordinary KMS positive functional and
$\pi_\o\cdot\g$ extends to an inner automorphism of 
$\cM\equiv\pi_\o(\gA)''$, implemented by a selfadjoint unitary 
$\G$ in the centralizer $\cM_{\o}$ of $\o$.
Moreover $\f$ is proportional to $\o(\G\cdot)$.
\end{Prop}
\begin{Cor}
If $\f$ is the unique non-zero graded KMS functional (up to a phase)
of $\gA$, then $\o = |\f|$ is extremal KMS, i.e. 
$\cM=\pi_{\o}(\gA)''$ is a factor.
\end{Cor}
\begin{proof}
As usual $\gA$ is identified with $\pi_\o(\gA)$.
If $Z(\cM)\neq \Complessi$, there exist two non-zero projections
$z_1,z_{2}\in Z(\cM)$ with sum $1$. Thus $\f(z_1\cdot),\f(z_2\cdot)$
are different graded KMS functionals on $\gA$. This can be checked 
since both the extensions of $\a$ and $\g$ act trivially on $Z(\cM)$, 
and by usual approximation arguments. By the uniqueness assumption 
there exists a constant $\l$ with $\f(z_1 a)=\l \f(z_2 a),\, a\in \gA$, 
then by continuity the same equality holds for $a\in \cM$, 
thus $\f(z_1 a)= \f(z_1 z_1  a) =\l \f(z_2 z_1 a) =0$. Analogously
$\f(z_2\cdot)=0$, thus $\f=0$.
\end{proof}
For our purposes Proposition \ref{BL} allows us to consider ordinary
KMS states instead of general graded KMS functionals.

Let $\r$ be  a graded endomorphism of $\gA$, namely an endomorphism of $\gA$ 
commuting with $\g$. 
In this graded context, we shall say that $\r$ is covariant if there exists a  covariance
$\a$-cocycle of unitaries $u(t)\in \gA$ such that $\g(u(t))=u(t)$.

In the following we again identify $\gA$ with 
its image $\pi_\o(\gA)$.

\begin{Lemma}\label{mod} Let $\f$ be a graded KMS functional for $\a$ and 
$\r$ a graded covariant endomorphism of $\gA$ as above. If $\cM$ is a 
factor and $\r$ extends to a finite index irreducible 
endomorphism of $\cM$,
then $\r$ has finite holomorphic dimension, both with respect to $\o$ 
and $\f$. Indeed, if $\o(u(\cdot))$ is continuous,
\[\label{eqholdim}
d_\f(u_\r )=\pm d_\o(u_\r ) \ ,
\]
namely $\f(1)d_\o(u) =\pm\acb\f(u(t))$.
\end{Lemma}
\begin{proof} The holomorphic property of $u$ is a direct consequence
of the holomorphic property of the Connes cocycle 
because $u$ is indeed a Connes Radon-Nikodym cocycle, up to phase, 
with respect to two bounded positive normal functionals 
of $\cM$ (see \cite{L4} and the previous section).

Note now that, since the extension of $\r$ to $\cM$ (still denoted by
$\r$) commutes with $\g={\rm Ad}\G$, we have $\r(\G)\G^*\in
\r(\cM)'\cap\cM =\mathbb C$, thus $\r(\G)=\pm \G$ because
$\G$ is self-adjoint, hence $\Phi_\r(\G)=\pm\G$, where
$\Phi_\r$ is the (unique) left inverse of $\r$.

To check eq. (\ref{eqholdim}), 
recall that, by the holomorphic properties of Connes cocycles,
we have (see \cite{L4}):
\begin{equation}
\acb \o(Xu(t))=d_\o(u)\o\cdot\Phi_\r(X), \quad \forall 
X\in\cM\ .
\end{equation}
As above we have
the polar decomposition $\f =\o(\Gamma\cdot)$. Then
\[
\acb \f(u(t))=\acb \o(\Gamma u(t)) =
d_\o(u)\o\cdot\Phi_\r(\G)=
\pm d_\o(u)\o(\G)=\pm \f(1)d_\o(u),
\]
namely the holomorphic dimension with respect to $\f$ coincides with 
the holomorphic dimension with respect to $\o$, up to a sign.
\end{proof}
Because of the above Lemma \ref{mod}, it is more convenient
to define the holomorphic dimension $d_\f(u)$ directly 
with respect to the modulus $\o$ of $\f$.

Before concluding this section, we recall that 
the interest in (bounded) graded KMS functionals
is limited by the following no-go theorem.

\begin{Prop}\emph{(\cite{BL})}\label{aa}
Let $\f$ be a graded KMS functional of $\gA$ with respect to $\a$
as above. If there exists
a $\f$-asym\-pto\-ti\-cally abelian sequence $\b_n\in{\rm Aut}(\gA)$, then
$\g=\iota$ and
${\f}$ is an ordinary KMS functional.
\end{Prop}
Here the $\b_n$'s commute with the grading and the 
$\f$-asym\-pto\-ti\-cally abelianness means that $\f\cdot\b_n=\f$
and $\f(c[\b_n(a),b])\to 0$ for all $a,b,c\in\gA$, where the 
commutator is a graded commutator.
\smallskip
\subsubsection{Table of dimensions.}
Before concluding this section we display the following table that
summarizes the various notions of dimension we are dealing with.
\bigskip

{\footnotesize \begin{tabular}{l l l}
\textit{Dimension}  &  \textit{Definition}&  \textit{Context}\\
\phantom{-}&\phantom{-}&\phantom{-}\\
Intrinsic   &  $d(\r) = \sqrt{||R_{\r}||\,||\bar R_{\r}||}$, 
(standard $R_{\r},\,\bar R_{\r}$)& Tensor C$^*$-categories\\
\phantom{-}&\phantom{-}&\phantom{-}\\
Analytical  &  $d_{an}(\r) = \sqrt{[\cM:\r(\cM)]}$& Subfactors\\
\phantom{-}&\phantom{-}&\phantom{-}\\  
Statistical &  $d_{DHR}(\r)=|\Phi_{\r}(\varepsilon_{\r})|^{-1}$,
($\varepsilon_{\r}$ stat. operator) & QFT, localized endomorphisms\\ 
\phantom{-}&\phantom{-}&\phantom{-}\\
Holomorphic & $d_\f(u)=
\underset{t\to i\b}{\textnormal{anal.cont.\,}}\f(u(t))$& Unitary 
cocycles\\
\phantom{-}&\phantom{-}&\phantom{-}\\
Geometric &  $d_{geo}(\r)=
\sqrt{d_\f(u_\r )d_\f(u_{\bar\r} )}$& Covariant endomorphisms\\           
\end{tabular}}
\bigskip

Here $[\cM:\r(\cM)]$ denotes the minimal index, namely the Jones
index with respect to the minimal expectation. 
We have omitted the notion of minimal dimension  
$d_{min}(\r)\equiv\textnormal{min}
\sqrt{||R_{\r}||\,||\bar R_{\r}||}$, in the context of C$^*$-tensor 
categories, as it turns out to coincide with the intrinsic 
dimension \cite{R}.

\section{The che\-mi\-cal po\-ten\-tial in Quan\-tum Field Theory.}
\label{QFTM}
In this section we examine certain aspects of the chemical potential 
for thermal states in Quantum Field Theory. Our discussion will rely
on basic results as the description of the chemical potential in terms
of extensions of KMS states \cite{AHKT} and
the construction of the field net and the gauge group \cite{DR}. 
Together with certain results for 
tensor categories \cite{LRo}, our 
analysis will show a splitting of the incremental free energy 
into an absolute part, that depends only on the charge and not on the 
state, and a part which is asymmetric with respect to the charge conjugation;
this indeed represents the chemical potential that labels the equilibrium states.

Let $\mathbb M$ be the Minkowski spacetime $\Reali^{d+1}$, with 
$d\geq 2$, and $\gA$ a net of local observable von Neumann 
algebras on $\mathbb M$, namely we have an 
inclusion preserving map
\[
\cO\to\gA(\cO)
\]
from the set $\cK$ of (open, non-empty) double cones 
of $\mathbb M$ to von Neumann algebras
$\gA(\cO)$ on a Hilbert space $\cH$. 
If $E\subset\mathbb M$ is arbitrary, we set $\gA(E)$ for the
C$^*$-algebra generated by the von Neumann algebras $\gA(\cO)$
as $\cO\in\cK$ varies  $\cO\subset E$
($\gA(E)=\mathbb C$ if the interior of $E$ is empty) and denote
by $\gA \equiv \gA(\mathbb M)$ the quasi-local C$^*$-algebra.
We denote the quasi-local C$^*$-algebra and the net itself by the same 
symbol, but this should not create confusion.
We shall also denote by $\cA(E)=\gA(E)''$ the von Neumann 
algebra generated by $\gA(E)$ (of course $\gA(\cO)=\cA(\cO)$ if 
$\cO\in\cK$).

We assume the following properties\footnote{These properties 
automatically hold, in particular, in a Wightman theory \cite{BW}. 
Here we will also consider the case of a reducible net.}:
\smallskip

{\it Additivity}:   If $\cO,\cO_1,\dots, \cO_n$ are 
double cones and $\cO_1\cup\dots\cup\cO_n\subset \cO$,
then $\gA(\cO_1)\vee\dots\vee\gA(\cO_n)\subset\gA(\cO)$.

Here and in the following, the lattice symbol $\vee$ denotes the von 
Neumann algebra generated.
\smallskip

{\it Wedge duality} (or essential duality):   
If $W\subset \mathbb M$ is a wedge region (namely a Poincar\'e 
translate of the region $\{x\in\mathbb M: x_1>x_0\}$), then
$$
\cA(W') = \cA(W)^c\ .
$$
Here $\cN^c$ denotes the relative commutant in the von Neumann algebra
$\cM\equiv\gA''$, namely $\cN^c\equiv \cN'\cap\cM$, and $W'$ denotes 
the spacelike complement of $W$.

In particular the net $\gA$ is {\it local}, namely $\gA(\cO_1)$
and $\gA(\cO_2)$ commute if the double cones $\cO_1$ and 
$\cO_2$ are space-like separated. 
As in the case of irreducible nets,
one may consider the dual net, here defined as
\[
\gA^d(\cO)\equiv \gA(\cO')^c, \quad\cO\in\cK\ ,
\]
and show the following, \cite{R2}:
\begin{Prop} $\gA^d$ satisfies {\it Haag duality}, by which we mean 
here that
\[
\gA^d(\cO) = \gA^d(\cO')^c, \quad\cO\in\cK\ ,
\]
where $\gA^d(\cO')$ is the C$^*$-algebra associated to $\cO'$ in
the net $\gA^d$. In particular $\gA^d$ is local.

Moreover $\gA$ and $\gA^d$ have the same weak closure:
\[\gA''=(\gA^d)''=\cM\ .
\]

\end{Prop}
\begin{proof} Because of additivity and wedge duality one can write
\begin{equation}\label{dual}
\gA^d(\cO)=\bigcap_{W\supset \cO}\cA(W)
\end{equation}
(intersection over all wedges containing $\cO$).

If $\cO_1\in\cK$ is spacelike separated from $\cO\in\cK$ there is a wedge $W$
with $W\supset\cO$ and  $W'\supset\cO_1$, hence $\gA^d$ is local.
As $\gA^d$ extends $\gA$ and is local, wedge duality must hold for
$\gA^d$ too, therefore $\gA^d(W)=\gA(W)$ for all wedges $W$.
In particular the global von Neumann algebra associated with
$\gA^d$ coincides with the one associated with $\gA$:
\[
(\bigvee_{\cO} \gA^d(\cO))''=(\bigvee_{W} \gA^d(W))''
=(\bigvee_{W} \gA(W))''=\cM .
\]
Analogously, it follows from formula \ref{dual} that $\gA(\cO')''=\gA^d(\cO')''$, 
namely $\gA^d(\cO)=\gA^d(\cO')^c$, 
that is to say Haag duality holds for $\gA^d$.
\end{proof}
Note however that, since $\cM$ is not a type I factor in general,
$\gA(\cO')''$ may be non-normal in $\cM$, namely
\[
\gA^d(\cO)^c = \gA(\cO')^{cc}
\]
may be strictly larger than $\gA(\cO')''$ if $\cO\in\cK$.
\smallskip

{\it Translation covariance}: There exists a
unitary representation $U$ of $\Reali^{d+1}$ making $\gA$ covariant:
\[
\t_x(\gA(\cO))=\gA(\cO + x),\ x\in\Reali^{d+1},\ \cO\in\cK\ ,
\]
where we have set $\t_x\equiv{\rm Ad}U(x).$
\smallskip

{\it Properly infiniteness and Borchers property B}
\footnote{In the vacuum representation these properties follows by 
positivity of the energy, see \cite{L,H}.}: If $\cO\in\cK$, then $\gA(\cO)$ is a 
properly infinite von Neumann algebra. If $\cO,\tilde\cO$ are double cones 
and $\cO+x\subset\tilde\cO$ for $x$ in 
a neighborhood of $0$ in $\mathbb R^{d+1}$, then every non-zero  projection
$E\in\gA(\cO)$ is equivalent to $1$ in $\gA(\tilde\cO)$. 
\smallskip

{\it Factoriality}: $\cM\equiv\gA''$ is a factor.
\smallskip

\noindent
A {\it localized endomorphism} $\r$ of $\gA$ 
is an endomorphism of $\gA$ 
such that $\r|_{\gA(\cO')}={\rm id}|_{\gA(\cO')}$
for some $\cO\subset\cK$. Two localized endomorphisms $\r,\r'$ 
are equivalent ($\r\simeq\r'$) if there is a unitary  $u\in\cM$ 
such that $\r'= {\rm Ad}u\cdot\r$. In the DHR theory \cite{DHR} 
$\cM= B(\cH)$, namely $\gA$ is irreducible, but most of what we are 
saying holds in the reducible case as well.

A localized endomorphism $\r$ is translation covariant if
there exists a $\tau$-cocycle of unitaries $u(\r,x)\in\cM$ such that
\begin{equation}\label{cov.formula}
{\rm Ad}u(\r,x)\cdot \r_x =\r , \quad x\in \mathbb R^{d+1},
\end{equation}
where $\r_x\equiv \t_x\cdot\r\cdot\t_{-x}$.

The equivalence classes of irreducible, translation covariant,
localized endomorphisms are the {\it superselection sectors} of $\gA$.

The translation covariant 
endomorphisms of $\gA$ form a tensor category.  
If $T\in (\r,\r')$ is an intertwiner between
$\r$, $\r'$, namely $T\in\cM$ and $T\r(X)=\r'(X)T$ for all $X\in\gA$,
then, as an immediate consequence of the Haag duality property for
$\gA^d$, we have
 $T\in\gA^d(\cO)$ if $\cO\in\cK$ contains the localization regions 
of both $\r$ and $\r'$. Now the unitary intertwiners changing the 
localization region of $\r$ (in particular the unitaries $u(\r,x)$ 
above) are used to define Roberts cohomology \cite{R}. In particular the
endomorphism $\r$ can be reconstructed from these charge transfers;
as they are local operators in $\gA^d$, they also provide an extension 
of $\r$ to $\gA^d$ with the same localization (only in the case 
$d\geq 2$). The superselection structure for $\gA$ and $\gA^d$
coincide (the extension map is a full functor) and, replacing $\gA$
by $\gA^d$, we may thus assume that $\gA$ satisfies Haag duality
(in the original representation of $\gA$).

As shown in \cite{DHR}, attached with any localized 
endomorphism $\r$ there is a 
unitary representation of the permutation
group $\mathbb P_{\infty}$, the statistics of $\r$, that is 
classified by a statistics parameter $\lambda_\r$ whose possible 
values are $\lambda_\r= 0,\pm 1,\pm \frac{1}{2},\pm \frac{1}{3}\dots$. 
Thus the {\it statistical 
dimension} $d_{DHR}(\r)=|\lambda_\r|^{-1}$ takes integral values
\[
d_{DHR}(\r)= 1,2,3,\dots +\infty\ .
\]
By the index-statistics theorem \cite{L1,L2}, $d_{DHR}(\r)$ 
coincides with an analytic 
dimension, the square root of the Jones index
\[
d_{DHR}(\r)=[\gA:\r(\gA)]^{\frac{1}{2}}\ ,
\]
(one way to read $[\gA:\r(\gA)]$ is $[\gA(W):\r(\gA(W))]$, with $W$ a 
wedge region).

We shall denote by $\cL$ be the tensor category of translation
covariant localized endomorphisms of $\gA$ with finite statistics 
(i.e. with finite dimension). For an object $\r$ of $\cL$, the 
intrinsic dimension coincides with the statistical dimension 
\cite{L1}:
\[
d(\r)=d_{DHR}(\r)\ .
\]
\begin{Thm}\label{hol4dim} Let  $\gA$ be as above,
$\a_t = \t_{x(t)}$ a one-parameter automorphism group of translations 
of $\gA$ and $\f$ a translation invariant state which is 
extremal KMS  for $\a$.

If $u(\r,t)$ is a $\a$-covariance cocycle for the irreducible localized 
endomorphism $\rho$ (i.e. eq. (\ref{cov}) holds), 
then $u_\r $ is holomorphic $	d(\r)\leq d_{geo}(\r)$.

If moreover $\pi_\f$ satisfies Haag duality, then
\begin{equation}
	d(\r)=d_{geo}(\r)\ .
\end{equation}
In particular the right hand side in the above formula is independent 
of the KMS state $\f$.
\end{Thm}
As is known the cyclic vector $\x_\f$ in the GNS representation of a 
KMS state $\f$ is separating for 
$\pi_\f(\gA)''$, namely $\f$ is a separating state. 
This is crucial for the following theorem of Takesaki
and Winnink.
\begin{Thm}\label{Takesaki}{\rm (\cite{TW})}. 
A KMS state $\f$ is locally normal, namely $\f|_{\gA(\cO)}$ is 
normal for any double cone $\cO\in\cK$.
\end{Thm}
Note that if $\gA$ is a Poincar\'e covariant net and $\r$ a 
Poincar\'e covariant localized endomorphism, then the covariance cocycle 
$u(\r,L) \ (\r\in\cL, L\in\Ppo)$ is uniquely fixed  because
$\Ppo$ has no non-trivial finite-dimensional unitary representation, 
hence it is a two-variable cocycle. In particular, if $\r$ extends to an irreducible endomorphism of the weak closure, then
$d(\r)=d_{geo}(\r)$ by Prop. \ref{2v}. We shall return on this point in the 
next section.

If moreover there is a PCT symmetry $j$ for $\gA$, or an 
anti-automorphism $j$ of $\gA$ such that $j^{-1}\cdot\r\cdot j=\bar\r$,
and $\f\cdot j=\f$, then $d(\r)=d_\f(u_\r )$ by Prop. \ref{PCT}.

We have essentially mentioned that the tensor category $\cL$ has a permutation 
symmetry (if dim($\mathbb M)\geq 3)$, as shown in \cite{DHR}. 
Then, by \cite{DR} there exists a field 
net $\gF$ of von Neumann algebras $\gF(\cO)\supset\gA(\cO)$, with 
normal commutation relations, with $\gA = \gF^G$ the fixed point
of $\gF$ under the action $\g$ of a compact group $G$ of internal 
symmetries of $\gF$. One has $\gA'\cap\gF=\mathbb C$, where 
$\gF$ is the quasi-local field C$^*$-algebra $\gF(\mathbb M)$.
Every endomorphism $\rho\in\cL$ 
is implemented by a Hilbert space of 
isometries in $\gF$.

We now relax the pointwise continuity condition in \cite{AHKT}.
We denote by $\tilde\t$ the translation automorphism group on $\gF$ 
extending $\t$ and by $\tilde\a=\tilde\t_{x(\cdot)}$ the one-parameter 
automorphism group extending $\a$.
\begin{Lemma} Let $\f$ be an extremal KMS state of $\gA$
with respect to $\a$.
There exists a locally normal state $\psi$ of $\gF$ that 
extends $\f$ and is extremal KMS with respect to
$\tilde\a_t\cdot \g_{g(t)}$, with $t\to g(t)$ a one-parameter subgroup
of $G$. 
\end{Lemma}
\begin{proof} 
Let $\gF_c\subset\gF$ denote the sub-C$^*$-algebra of all elements with
pointwise norm continuous orbit under the action of $\tilde\t\cdot
\g$ of $\Reali^{d+1}\times G$, and set 
$\gF_c(\cO)\equiv \gF_c\cap\gF(\cO)$. 
We have
\[
{\gF_c(\cO)}^{-} \supset\gF(\cO_0),\quad \forall\ \cO_0\in\cK,\
\bar\cO_0\subset\cO,
\]
where ${\gF_c(\cO)}^{-}$ is the $\s$-weak closure of $\gF(\cO)$. Indeed 
if $X\in\gF(\cO_0)$ and $j_n$ is an approximation of the identity in 
$\Reali^{d+1}$ by continuous functions with support in a ball of radius
$\frac{1}{n}$, then
\[
X_n\equiv\int j_n(X)\tilde\t_x(X){\rm d}x\to X
\]
($\s$-weak convergence) and $X_n$ has pointwise $\tilde\t$-orbit and,
for large $n$, belongs to $\cF(\cO)$ (the $\g$-continuity is checked
similarly).

Set $\f_c\equiv \f|_{\gA_c}$, where $\gA_c\equiv\gF_c\cap\gA$
and $\tilde\f=\f\cdot\varepsilon$, where $\varepsilon=
\int_G \g_g{\rm d}g$ is the expectation of $\gF$ onto $\gA$ as above.
Clearly $\tilde\f$ is a $\tilde\t$-invariant locally normal state
of $\gF$ and so is its restriction $\tilde\f_c$ to $\gF_c$.
Let $\psi_c\prec\tilde\f_c$ be an extremal $\tilde\a$-invariant 
state of $\gF_c$ extending $\f_c$.
By the AHKT theorem \cite{AHKT} 
$\psi_c$ is KMS with respect to a one-parameter automorphism group 
$t\to \tilde\a_t\cdot \g_{g(t)}$ of $\gF_c$.
Since $\tilde\f_c$ is locally normal and dominates $\psi_c$, 
also $\psi_c$ is locally normal, thus it extends to a locally normal
state $\psi$ of $\gF$. By usual arguments, $\psi$ is a KMS
state on $\gF$ with respect to $\tilde\a$.
\end{proof}
\begin{Lemma}\emph{(\cite{AHKT} \textrm{Prop. III.3.2})}\label{separ}
Let $\gA\subset\gF$ 
be C$^*$-algebras, $\f$ a state of $\gA$ and $\tilde\f$ an extension 
to $\gF$. If $\tilde\f$ is separating, then $\pi_{\tilde\f}|_{\gA}$
is quasi-equivalent to $\pi_\f$.
\end{Lemma}
\begin{Cor}\label{AHKT} Let $\gF$ be a C$^*$-algebra and 
$\r$ an inner endomorphism of $\gF$. If $\f$ is a separating 
state of $\gF$, the GNS representations $\pi_\f$ 
and $\pi_{\f\cdot\r}$ are quasi-equivalent.
\end{Cor}
\begin{proof}
Let $H\subset\gF$ be a Hilbert space of isometries implementing $\r$ 
on $\gF$ and let $\{v_i,\, i=1,\dots,n\}$ be an orthonormal basis 
of $H$, thus the 
$v_i$'s are isometries in $\gF$ with orthogonal final projections 
summing up to the identity and $\r(X)=\sum_i v_i X v_i^*$, $X\in\gF$. 
Then
\begin{multline*}
(\pi_{\f\cdot\r}(X)\xi_{\f\cdot\r},\xi_{\f\cdot\r})=\f(\r(X))=
(\pi_\f(\r(X))\xi_\f,\xi_\f)\\
=\sum(\pi_{\f}(X)\xi_i,\xi_i)=(\pi_{\f}(X)\oplus\dots
\oplus\pi_{\f}(X)\bar\xi,\bar\xi),\quad X\in\gF,
\end{multline*}
where $\xi_i=\pi_{\f}(v_i^*)\xi_\f$ and $\bar\xi=\oplus\xi_i$. Therefore
$\pi_{\f\cdot\r}\prec\pi_{\f}\oplus\dots\oplus
\pi_{\f}$. On the other hand $\bar\xi$ is a separating 
vector for $\pi_\f\oplus\dots\oplus\pi_\f(\gF)''$; indeed elements of
$\pi_\f\oplus\dots\oplus\pi_\f(\gF)''$ have the form 
$Y=X\oplus\dots\oplus X$, $X\in\pi_\f(\gF)''$, thus $Y\bar\xi=0$ iff
$X\pi_{\f}(v_i^*)\xi_\f=X\x_i=0,\ \forall i$, thus iff 
$X\pi_{\f}(v_i^*)=0$ because $\xi_\f$ is separating, which is 
equivalent to $X=0$ because $\sum v_i^*v_i =0$.
\end{proof}
\begin{Lemma}\label{sep} 
Let $\f$ be an extremal KMS state of $\gA$.
The GNS representations 
$\pi_\f$ and $\pi_{\f\cdot\r}$ are quasi-equivalent.
\end{Lemma}
\begin{proof}
Let  $\psi$ a KMS state of $\gF$ extending $\f$. As $\psi$
is a separating state of $\gF$, also $\psi\cdot\r_H$ is also separating 
as in the proof of Prop. \ref{AHKT}, where $\r_H$ is the inner 
endomorphism of $\gF$ implemented by $H$. 
As $\psi\cdot\r_H$ extends $\f\cdot\r$ we have 
by Lemma \ref{AHKT} that $\pi_\f\simeq\pi_{\psi}|_{\gA}
\simeq\pi_{\psi\r_H}|_{\gA}\simeq\pi_{\f\cdot\r}$, where the symbol
``$\simeq$'' denotes quasi-equivalence.
\end{proof}
\begin{Cor} \label{ext4dim}
Every localized endomorphism $\rho\in\cL$ is normal with respect to
$\f$.
\end{Cor}
\begin{proof}
The proof now follows by Lemma \ref{sep} and Lemma \ref{ext}.
\end{proof}
{\it Proof of Theorem} \ref{hol4dim} As $\f$ is extremal KMS, 
$\cM\equiv\pi_\f(\gA)''$ is a factor. If $\rho$ is a localized endomorphism,
by Corollary \ref{ext4dim} both $\rho$ and its conjugate extend to
$\cM$. The conjugate equation then holds on $\cM$ showing that 
the extension $\hat\r$ of $\r$ to $\cM$ has finite dimension. 
We have $d_{an}(\hat\r)=d(\r)$
(a priori we only have $d_{an}(\hat\r)\leq d(\r)$). Indeed 
Proposition \ref{extbraided} applies. 

If Haag duality holds in the representation $\pi_\f$, then 
by the following Lemma \ref{full} $\hat\r$ is irreducible if 
$\r$ is irreducible, therefore the last part of the statement follows
by Prop. \ref{dim}.
The rest now follows from the analysis in the previous section.
\qed
\begin{Lemma}\label{full}
With the above notation, if $\pi_\f$ satisfies Haag duality, then 
the extension map $\r\in\cL\to\hat\r\in{\rm End}(\cM)$ is a full 
functor.
\end{Lemma}
\begin{proof} With $\r\in \cL$  irreducible, we have to show that
$\hat\r$ is irreducible too. Let $T\in(\hat\r,\hat\r)$, namely
$T\in\cM$ and $T\hat\r(X)=\hat\r(X)T$ for all $X$ in $\cM$.
As $\r$ is localized in a double cone $\cO$, $\r$ acts identically on 
$\gA(\cO')$, hence $T\in \gA(\cO')'\cap\cM=\gA(\cO)$,
thus $T\in(\r,\r)=\mathbb C$ as desired.
\end{proof}
\subsection{The absolute and the relative part of the incremental 
free energy.}
\label{abschem}
Beside the description of the chemical potential in terms of extensions
of KMS states, the AHKT work provides an intrinsic description
of the chemical potential within the observable algebra \cite{AHKT}, see also \cite{H}, that 
was made explicit only in the case of abelian charges 
(automorphisms).

Let's recall this point. Let 
$\gA$ be a unital C$^*$-algebra with trivial centre and $\a$ a one-parameter automorphism group 
as before and $\r$ a covariant automorphism of $\gA$, thus 
$\r\cdot\a_t\cdot\r^{-1}=\textnormal{Ad}u(t)\cdot\a_t$ for some 
$\a$-cocycle of unitaries $u(t)\in\gA$. Notice that $u$ is unique up 
to multiplication by a one dimensional character of $\Reali$, that one
fixes once for all.

If now $\f$ is an extremal KMS state for $\a$ such that 
$\f\cdot\r^{-1}$ 
and $\f$ are quasi-equivalent, thus $\r$ extends to the factor 
$\cM\equiv \pi_\f(\gA)''$, then 
\begin{equation}
	u(t)= e^{-i\m_{\r}(\f)t}(D\f\cdot\r^{-1}:D\f)_{-\b^{-1}t}
	\label{chem}
\end{equation}
for some $\m_{\r}(\f)\in\Reali$, called the chemical potential of $\f$
(we are assuming that $\f(u(\cdot))$ is continuous). 
A relevant observation is that, although
$\m_{\r}(\f)$ depends on the initial phase fixing for $u$, 
\[
\m_{\r}(\f'|\f)\equiv\m_{\r}(\f') -\m_{\r}(\f)
\]
is independent of that and is therefore an intrinsic quantity 
associated with a pair $\f$, $\f'$ of extremal 
KMS states. In other words the chemical potential is a label for the 
different extremal KMS states.

Now the above argument goes true in more 
generality if $\r$ is an endomorphism of $\gA$ that extends to a 
finite-index irreducible endomorphism $\cM$, once we replace
$\f\cdot\r^{-1}$ with $\f_\r\equiv\f\cdot\Phi_\r$, 
where $\Phi_\r$ is the minimal left inverse of the extension of $\r$.
In general, for a given charge $\r$, the chemical 
potential is only defined with respect to the two thermal states
\[
\mu_\r(\f'|\f)\equiv\b^{-1}\log d_\f(u) -\b^{-1}\log d_{\f'}(u)=
 \b^{-1}\log d_{\f'}((D\f_\r:D\f))\ .
\]
Moreover, if there is a canonical way to choose the cocycle $u$, independently
of the state $\f$, then 
\begin{equation}\label{absoche}
\m_\r(\f)\equiv \b^{-1}\log d_\f (u) - \b^{-1}\log d(u)\ .
\end{equation}
defines an {\it absolute}
chemical potential in the state $\f$, associated with the charge $\r$. 

The above discussion relies of course on the normality of the 
endomorphism $\r$ in a thermal state, a deep fact, proved 
in \cite{AHKT} when the endomorphisms $\r$ are 
associated to the dual of a compact gauge group, with certain
asympotically abelian and cluster properties for the dynamics.

We now apply the above discussion to the case of quantum
relativistic statistical mechanics, namely we consider 
thermal states for the time evolution in a 
quantum field theory on Minkowski spacetime. In this situation,
there is a two variable cocycle for the Poincar\'e  covariant endomorphisms 
with finite dimension, which is unique because the Poincar\'e has no 
non trivial finite dimensional unitary representation.
Hence $\m_{\r}(\f)$ can be defined intrinsically by eq. (\ref{absoche}). 

To be more explicit, let $\gA$ be a net of von Neumann algebras 
on the Minkowski space, as 
in the previous section. We further assume that $\gA$ is 
Poincar\'e covariant, namely there is a unitary representation $U$
of $\Ppo$ on $\cH$ that acts covariantly on $\gA$ and extends the 
translation unitary group.

Our endomorphisms $\r$ are assumed to be covariant with respect to the action 
of $\Ppo$, namely there exists an $\a$ cocycle of unitaries 
$u(\r,L)\in\gA$ such that
\begin{equation}
\textnormal{Ad}u(\r,L)\cdot\r=\a_L\cdot\r\cdot\a^{-1}_{L}, 
\quad L\in\Ppo .
	\label{unicoc}
\end{equation}
where $\a_L={\rm Ad}U(L)$. This is indeed a two-variable cocycle.
 
Restricting this cocycle to the subgroup of time translation, 
we obtain a canonical choice for the unitary cocycle 
$u(\r,t), \, t\in\Reali$, for the one parameter group
$\a$, which is indeed a two-variable cocycle in $\cL\times\Reali$.

\begin{Thm}\label{decomposition}
Let $\gA$ be the quasi-local observable C$^*$-algebra and $\a$ 
a one-parameter (time) translation automorphism group as in 
the previous section.
If $\r$ is a 
Poincar\'e covariant irreducible localized 
endomorphism with finite dimension, we have:
\begin{itemize}
\item [$(i)$] If $\f$ is an extremal KMS state for $\a$ satisfying 
Haag duality,
there exists a
chemical potential $\mu_{\r}(\f)$ associated with $\r$ defined by 
the canonical splitting
\[
\log d_\f(u_\r )= \log d(\r) + \b\mu_{\r}(\f)
\]
and satisfies
\[
\mu_{\r}(\f)=-\mu_{\bar\r}(\f)\ .
\]
The intrinsic dimension is thus given by
\[
\log d(\r) =\frac{1}{2}(\log d_\f(u_\r )+
\log d_\f(u_{\bar\r}))
\]
independently of $\f$.
\item [$(ii)$] If there is 
a time reversal symmetry as in Prop. \ref{PCT} and $\f\cdot j =\f$, then
$\mu_\r(\f)=0$, namely
$$
d_\f(u_\r )=\log d(\r),
$$
independently of $\f$.
\end{itemize}
\end{Thm}
\begin{proof} As above noticed, $u(\r,L)$ is a two-variable cocycle 
for the action of $\cL_0\times\Ppo$ on $\gA$, where $\cL_0\subset\cL$ 
is the tensor category with conjugates generated by $\r$ (see 
\cite{L4}). Hence by \cite{L4}
\[
u(\r,L)^{\bullet}=u(\bar\r,L)\ .
\]
On the other hand, by the results in the previous section,  we may 
extend $\r$ to the weak closure of $\gA$ in the GNS representation 
of $\f$, and if then compare with the Connes cocycle, we have
\[
u(\r,t)=e^{i\mu_\r(\f)t}d(\r)^{-i\b^{-1}t}(D\f_\r:D\f)_{-\b^{-1}t} \ .
\]
But 
$(d(\r)^{it}(D\f_\r:D\f)_t)^{\bullet}=(d(\r)^{it}(D\f_{\bar\r}:D\f)_t)$,
hence
\[
u(\bar\r,t)=
e^{-i\mu_\r(\f)t}d(\r)^{-i\b^{-1}t}(D\f_{\bar\r}:D\f)_{-\b^{-1}t} \,
\]
namely $\mu_{\bar\r}(\f)=-\mu_\r(\f)$.

The last point is a consequence of Proposition \ref{PCT}.
\end{proof}
Note in particular that by the point $(ii)$ in the above theorem the 
chemical potential vanishes if there is a time-reversal symmetry,
therefore a non-trivial chemical potential sets an arrow of time,
in accordance with the second principle of thermodynamics.

Now we make contact with the analysis in \cite{L4}. Since the 
covariance cocycle
$u=u_\r $ is canonically defined in the above context,
once we choose the thermal state $\f$, the Hamiltonian
in the state $\f_\r$ is canonically defined as
\[
H_\r\equiv -i\frac{\rm d}{{\rm d}t}u(\r,t){e^{itH}}|_{t=0}\ ,
\]
where $H=\b^{-1}\log \Delta_{\xi}$  is the Hamiltonian
in state $\f$. Here $\xi$ is the GNS vector associated with $\f$.
The increment of the free energy between the states $\f$ and $\f_\r$
is then defined (cf. \cite{L4}) as
\[
F(\f|\f_\r)\equiv \f_\r(H_\r)-\b^{-1}S(\f|\f_\r)=-\b^{-1}\log(e^{-\b 
H_\r}\xi,\xi)\ ,
\]
where $S(\f|\f_\r)$ is Araki's relative entropy.
It is immediate from the last expression that
\[
F(\f|\f_\r)=-\b^{-1}\log d_\f(u_\r )=-\b^{-1}\log d(\r)-
\mu_\r(\f)\ .
\]
If $\f'$ is another extremal KMS state as above, so that $\r$ is 
normal with respect to $\f'$, we may now define
the increment of the free energy 
$F(\f'|\f'_\r)=-\b^{-1}\log d_{\f'}(u_\r )$,
where $\f'_\r=\f'\cdot\Phi_\r$, so we have:
\[
F(\f'|\f'_\r)-F(\f|\f_\r)=\mu_\r(\f) - \mu_\r(\f')=\mu(\f'|\f)\ ,
\]
moreover
\[
S_c(\r)=\log d(\r)^2=-\b(F(\f|\f_\r)+F(\f|\f_{\bar\r}))
\]
is an integer independent of $\f$. Here $S_c(\r)$ is the conditional 
entropy of $\r$, see \cite{L4}.

According to the thermodynamical formula ``${\rm d}F
={\rm d}E - T{\rm d}S$'', we have obtained the following relation:
\[
F(\f|\f_\r)=\mu_\r(\f) -\frac{1}{2}\b^{-1}S_c(\r) \ .
\]
The quantity $\mu_\r(\f)$ may be interpreted as part of the 
energy increment obtained by adding the the charge $\r$ to the 
identical charge, more specifically the part which is asymmetric with 
respect to charge conjugation. The total increment of the free energy
contains also a part which is symmetric under charge conjugation and
independent of the thermal equilibrium state, namely the intrinsic 
increment of the entropy $\frac{1}{2}S_c(\r)$ multiplied by the
temperature $\b^{-1}$.

The above analysis simply goes through when we consider the increment 
of the free energy between two thermal states $\f_\r$ and $\f_\s$ (cf. 
\cite{L4}). We shall make this explicitly in the context of Section 
\ref{blackhole}.
\section{Chemical potential. Low dimensional case.}
\label{CLD}
We now study the chemical potential structure in the low dimensional 
case. 
The higher dimensional methods \cite{AHKT}  cannot be
applied to this context, but we shall see that an analysis is possible
by using Wiesbrock's characterizations of conformal nets
on $S^1$ \cite{W}.
\subsection{KMS states and the generation of conformal nets.}
Let $\cI$ denote the set of all bounded open non-empty intervals of 
$\mathbb R$. We shall consider a net $\gA$ of von Neumann algebras 
on $\mathbb R$, namely an inclusion preserving map 
\begin{displaymath}
	I\in\cI\to\gA(I)
\end{displaymath}
from $\cI$ to von Neumann algebras $\gA(I)$, not necessarily acting
on the same Hilbert space.
We denote by the same symbol the 
quasi-local observable C$^{*}$-algebra 
$\gA=\cup_{I\in\cI}\gA(I)^{-}$  (norm closure).

We shall assume the following properties of $\gA$:
\smallskip

{\it a) Translation covariance:} There exists a one-parameter 
automorphism group $\t$ of $\gA$ that corresponds to the translations 
on $\mathbb R$,
\[
\t_s(\gA(I))=\gA(I+s),\quad I\in\cI, \, s\in\mathbb R\ .
\]

{\it b) Properly infiniteness:} For each $I\in\cI$, the von Neumann 
algebra $\gA(I)$ is properly infinite.\footnote{This assumption 
is needed only for the the local normality of
the KMS states (above Th. \ref{Takesaki} from \cite{TW}) and 
can alternatively be replaced by the factoriality of $\gA(I)$.}
\smallskip

Let now $\f$ be a KMS state on $\gA$ with respect to $\t$; for 
simplicity we set $\b=1$. Let $(\cH_\f,\pi_\f,\x_\f)$ be the 
associated GNS triple and $V$ the one-parameter unitary group 
implementing $\t$:
\[
V(s)\pi_\f(a)\x_\f=\pi_\f(\t_s(a))\x_\f\ .
\]
Note that by the KMS condition $s\to\f(a\t_s(b))$ is a continuous map 
for all $a,b\in\gA$, hence $V$ is strongly continuous.

Recall now that $\gA$ is {\it additive} (resp. {\it strongly additive}) if
\[
 \gA(I)\subset\gA(I_1)\vee\gA(I_2),
 \label{additivity}\]   
whenever $I, I_1, I_2 \in\cI$ and $I\subset I_1\cup I_2$ 
(resp. $I\subset \bar I_1\cup \bar I_2$), where the bar denotes the 
closure.

We now set $\cA(I)=\pi_\f(\gA(I))$, $I\in\cI$, which is a von Neumann 
algebra by Th. \ref{Takesaki}, and 
$\cA(E)\equiv\vee\{\cA(I): I\in\cI, I\subset E\}$ for any set 
$E\subset\mathbb R$ (the von Neumann algebras generated). 
Again we now assume $\pi_\f$ to be one-to-one
and identify $\gA(I)$ with $\cA(I)$, namely we consider the net 
already in its GNS representation.

The following KMS version of the Reeh-Schlieder theorem
is known to experts.

\begin{Prop} $\x$ is cyclic and separating for $\cA(I)$, 
if $I$ is a half-line. If $\gA$ is additive, then $\x$ is cyclic and separating
for all $\cA(I)$, $I\in\cI$.
\end{Prop} 
\begin{proof} Assume first that $I$ is a half-line and 
let $\eta\in\cH$ be orthogonal to 
$\cA(I)\x$; we have to show that $\eta=0$. Indeed if $I_0$ is a 
half-line and 
$\bar I_0\subset I$, then for all $a\in\cA(I_0)$
\[
 (\eta,V(s)a\x)=0,
 \]
 for all $s\in\mathbb R$ such that $I_0 + s\subset I$.
But, because of the KMS property, the function $s\to(\eta,V(s)a\x)$ is 
the boundary value of a function analytic in the strip 
$0<\textnormal{Im}z<\frac{1}{2}$ (as 
$V(\frac{i}{2})=\Delta^{\frac{1}{2}}$ and 
Dom$(\Delta^{\frac{1}{2}}\supset\cA(I_0)\xi$), 
hence it must vanish everywhere. It follows that
$\eta$ is orthogonal to $\cA(I_0 + s)\x$ for all $s\in\mathbb R$, 
hence $\eta$ is orthogonal to 
$\overline{\cup_{s\in\mathbb R}\cA(I_0 
+s)\x}\supset \overline{\gA\x}=\cH$.

Assume now that $I\in\cI$ and $\cA$ is additive. 
Set $\cA_0(I)=\vee\{\cA(I_0): I_0\in\cI, \bar I_0\subset I\}$.
We shall show that $\x$ is cyclic
for $\cA_0(I)$, hence for $\cA(I)$.

By the same argument as above 
$$
\eta\ \bot \
\cA_0(I)\x \Rightarrow \eta \ \bot \ \cA(I_0 + s)\x ,
\ \forall s\in\mathbb R, \ \bar I_0 \subset I,
$$ 
namely the orthogonal projection $P$ onto  $\overline{\cA_0(I+s)\x}$ is independent 
of $s\in\mathbb R$ and thus belongs to $\cap_s \cA_0(I+s)'=(\vee_s 
\cA_0(I+s))' = \gA'$ (by additivity). As $\x$ is separating for $\gA'$ 
and $P\x=\x$ it follows that $P=1$, namely $\overline{\cA_0(I)\x }= \cH$. 
\end{proof}

Now $\t$ extends to the rescaled modular group of 
$\cM\equiv\cA(\mathbb R)$ with respect 
to $\f$ and $\t_s(\cN)=\cA(s,\infty)\subset\cN,\ 
s>0$, where $\cN\equiv\cM(0,\infty)$, namely $(\cN\subset\cM,\x)$ is a 
half-sided modular inclusion of von Neumann algebras and by 
Wiesbrock's theorem \cite{W} 
there exists a $\x$-fixing one-parameter unitary group $U$ on $\cH$
with positive generator such that
\begin{gather}
V(s)U(t)V(-s)=U(e^s t)\\
U(1)\cM U(-1)=\cN
\end{gather}
Setting $\cB(a,b)\equiv\cA(\log a,\log b),\ b>a>0$, we have a net 
$\cB$ on the intervals of $(0,\infty)$ whose closure is contained in
$(0,\infty)$. We have the following, compare with \cite{BY}.
\begin{Prop}\label{tc}
 Let $\gA$ be an additive net as above and $\f$ a KMS state.
There exists a net $\cB$ on the intervals of $(0,\infty)$
such that $\cB(a,b)=\pi_\f(\gA(\log a,\log b))$ if $b>a>0$. 
$\cB$ is dilation covariant and $V$ is the dilation one parameter group.
$\cB$ is also translation  covariant with positive energy on 
half-lines,
namely there is a one-parameter $\x$-fixing unitary group $U$ with 
positive generator such that $\Ad U(t)\cB(a,\infty)= \cB(a+t,\infty)$,
where $\cB(a,\infty)=\vee_{b>a}\cB(a,b)$. 
\end{Prop}
\begin{proof} Clearly $V(s)\cB(a,b)V(-s)=\cB(e^s a, e^s b)$ for 
positive $a,b$. Setting $\bar\cB(a,\infty)\equiv\Ad U(a)\cM$
we have
\[
\Ad U(t)\bar\cB(a,\infty)=\bar\cB(t+a,\infty)\ ,
\]
therefore, by using the relation $V(s)U(t)V(-s)=U(e^st)$, it follows that
\[
\Ad V(s)\bar\cB(a,\infty)=\bar\cB(e^s a,\infty).
\]
On the other hand $\bar\cB(1,\infty)=\Ad U(1)\cM=\cN=\cB(1,\infty)$
hence
\[
\cB(e^s,\infty)=\Ad V(s)\cB(1,\infty)=
\Ad V(s)U(1)\cM=\Ad U(e^s)\cM=\bar\cB(e^s,\infty)\ ,
\]
showing the last part of the statement.
\end{proof}

We shall call the net $\cB$ the {\it thermal completion} of $\gA$ with
respect to $\f$. Note that the translation unitary group $V$ for $\cA$ becomes the 
dilation unitary group for $\cB$. 

Proposition \ref{tc} does not give the translation 
covariance of the net $\cB$ on the bounded intervals ($\cB(a,b)$ is 
not even defined if $a<0$). 

Further insight in the structure of the thermal completion net may be 
obtained by considering a {\it local} net $\gA$, namely
assuming the locality condition
\[
[\gA(I_1),\gA(I_2)]=\{0\}\ \textrm{if}\ I_1\cap I_2=\emptyset \ .
\]
 To construct a translation covariant net we 
define the following von Neumann algebras:
\begin{gather}\label{mtc}
\tilde\cB(0,1)=\bigvee_{s\leq 0}\Ad V_1(s)\cB(0,1),\\
\tilde\cB(0,a)=\Ad V(a)\tilde\cB(0,1),\ a\in\mathbb R,\\
\tilde\cB(a,b)=\Ad U(a)\tilde\cB(0,b-a), \ a,b\in\mathbb R\ .
\end{gather}
Here we have set $V_1(s)\equiv U(1)V(s)U(-1)$, the one-parameter
unitary group associated with the dilations with respect to  the point 
$1\in\Reali$.

From now on the net $\gA$ will be assumed to be local.

\begin{Thm}\label{ctc}
Let the net $\gA$ on the intervals of $\mathbb R$ be translation 
covariant, local, and additive and $\f$ a KMS state. 
With the above notations,  $\tilde\cB$ defines  a conformal net, 
indeed $\tilde \cB$ has a conformal extension to $S^1$.

As a consequence the dual net of $\tilde\cB$ on $\mathbb R$ is 
strongly additive and conformal.
\end{Thm}
We shall see that
\[
\tilde\cB(a,\infty)=\cB(a,\infty),\quad a\in\mathbb R,
\]
hence $\tilde\cB$ is an extension of $\cB$ on the positive 
half-lines and is conformal. We shall call $\tilde\cB$ the {\it
conformal thermal completion} of $\gA$ with respect to  $\f$.
\begin{proof}
We first show that the triple 
$\{\cB(0,\infty)',\tilde\cB(0,1),\cB(1,\infty),\x\}$ is a +hsm factorization
with respect to $\x$ in the sense of \cite{GLW}, 
namely these three algebras mutually commute
and $(\tilde\cB(0,1)\subset\cB(0,\infty),\x)$, $(\cB(1,\infty)\subset
\tilde\cB(0,1)',\x)$ and $(\cB(0,\infty)'\subset\cB(1,\infty)',\x)$
are +half-sided modular inclusions.

Now $\cB(0,\infty)'$ and $\cB(1,\infty)$ commute by the isotony of 
$\cA$; for the same reason $\cB(1,\infty)$ commute with
$\tilde\cB(0,1)$, indeed $\cB(1,\infty)$ is $\Ad V_1$-invariant
where, as above, $V_1(s)=U(1)V(s)U(-1)$.
Again $\cB(0,\infty)\supset\cB(0,1)$, hence $\cB(0,\infty)\supset\tilde
\cB(0,1)$ because $V_1(s)\cB(0,\infty)V_1(-s)\supset\cB(0,\infty)$
if $s\leq 0$ by translation-dilation covariance of $\cB$ on positive 
half-lines (Prop. \ref{tc}).

Concerning the hsm properties, the only non-trivial verification is that
$(\tilde\cB(0,1)\subset\cB(0,\infty),\x)$ is a +hsm inclusion, 
namely that
\[
\Ad V(s) \tilde\cB(0,1)\subset \tilde\cB(0,1),\ s<0.
\]
We thus need to show that
 for any fixed $t<0$ we have
\[
\Ad V(s) V_1(t)\cB(0,1)\subset \tilde\cB(0,1),\ s<0.
\]
Indeed if $s<0$ and $t<0$, there exist $s'<0$ and $t'<0$ such that
$V(s) V_1(t) = V_1(t')V(s')$, as follows immediately by the 
corresponding relation in the ``$ax+b$'' group. Therefore
\[
\Ad V(s) V_1(t)\cB(0,1)=
\Ad V_1(t')V(s')\cB(0,1)\subset\Ad V_1(t')\cB(0,1)\subset \tilde\cB(0,1)
\]
as desired.

By a result in \cite{GLW} there exists a conformal net $\tilde\cB$ on 
$\mathbb R$ such that the local von Neumann algebras associated 
 to $(-\infty,0)$, $(0,1)$ and $(1,\infty)$ are respectively
$\cB(0,\infty)'$, $\tilde\cB(0,1)$ and $\cB(1,\infty)$ and having $U$
and $V$ as translation and dilation unitary groups. 

By translation-dilation covariance, 
$\tilde\cB$ is then conformal thermal completion of $\cA$.
\end{proof}
We may also directly define the dual net $\cB^d$ of $\tilde\cB$ 
as the one associated 
with the half-sided modular factorization 
$(\cB(0,\infty)', \cB(1,\infty)'\cap\cB(0,\infty),
\cB(1,\infty),\x)$. This net is conformal, strongly additive and
\[
\cB^d(a,b)=\cB(a,\infty)\cap\cB(b,\infty)'.
\]
This is due to the equivalence between strong additivity and Haag 
duality on the real line for a conformal net, see \cite {GLW}.
Clearly we have
\[
\cB(a,b)\subset\tilde\cB(a,b)\subset\cB^d(a,b).
\]
thus $\cB^d$ is the dual net of $\tilde\cB$ and
\[
\cB^d=\tilde\cB\Leftrightarrow \tilde\cB \textnormal{ is strongly 
additive,}
\]
\begin{Cor}\label{commutant}
If $\cA$ is strongly additive
\[
\cB(1,\infty)'\cap\cB(0,\infty)=\bigvee_{s<0}V(s)\cB(0,1)V(-s)\ .
\]
\end{Cor}
\begin{proof} If $\cA$ is strongly additive, then $\cB$ is strongly 
additive (on the intervals of $(0,\infty)$), hence $\tilde\cB$ is 
strongly additive and the above comment applies.
\end{proof}
More directly, Corollary \ref{commutant} states that the relative 
commutant $\cA(0,\infty)'\cap\cM$ is the smallest von 
Neumann algebra containing $\cA(-\infty,0)$ which is mapped 
into itself by $\Ad\Delta^{it}$, $t>0$, where $\Delta$ is the modular 
operator associated with $(\cA(0,\infty),\x)$.

We shall say that the state $\f$ of $\gA$ satisfies {\it essential
duality} if
\[
\cA(0,\infty)'\cap\cM=\cA(-\infty,0)\ .
\]
We have:
\begin{Prop}
The following are equivalent:
\begin{itemize}
 \item[$(i)$] $\f$ satisfies essential duality,
\item[$(ii)$] For some (hence for all) $0<a<b$ we have 
$\cA(b,\infty)'\cap\cA(a,\infty)=\cA(a,b)$,
\item[$(iii)$]  $\cA$ is strongly additive and 
$\Ad\Delta^{it}\cA(-\infty,0)\subset\cA(-\infty,0)$ 
for all $t>0$, where $\Delta$ is the modular operator
of $(\cA(0,\infty),\x)$.
\end{itemize}
\end{Prop}
\begin{proof} $(i)\Leftrightarrow (iii)$ follows by the above comments.
On the other hand $(i)\Leftrightarrow (ii)$ because they are  
equivalent to the relative commutant property 
$\cB(a,b)=\cB(b,\infty)'\cap\cB(a,\infty)$ for $b>a$ and either 
$a=0$ or $a>0$, which are indeed equivalent conditions 
in the conformal case \cite{GLW}.
\end{proof}
\begin{Cor}
If $\f$ satisfies essential duality, then $\f$ satisfies Haag duality,
namely
\[
\cA(a,b)=(\cA(-\infty,a)\vee\cA(b,\infty))^c\ ,\quad a<b \ ,
\]
where $\cdot ^c$ denotes the relative commutant in $\cM$.
\end{Cor}
\begin{proof}
If $\f$ satisfies essential duality then, since 
$\cA(b,\infty)'\cap\cA(a,\infty)\subset\cM$ by $(ii)$ of the above
proposition, we have 
\[
\cA(a,b)=\cA(b,\infty)'\cap\cA(a,\infty)=
\cA(b,\infty)^c\cap\cA(a,\infty)
\]
for $b>a>0$. On the other hand, by essential duality, we have
$\cA(a,\infty)=\cA(-\infty,a)^c$, hence
\[
\cA(a,b)=\cA(-\infty,a)^c\cap\cA(b,\infty)^c=
(\cA(-\infty,a)\vee\cA(b,\infty))^c
\]
as desired. The case of arbitrary $a<b$ is obtained by translation 
covariance.
\end{proof}
Hence essential duality in a thermal state can occur only if the 
original net is strongly additive. It would interesting to see if the 
converse holds true, namely if all KMS states on a strongly additive net
satisfy essential duality.

Note also the, in contrast to the situation occurring in the vacuum 
representation, the equality
$\cA(a,b)'\cap\cM=\cA(-\infty,a)\vee\cA(b,\infty)$ cannot hold
in any thermal state, unless the superselection structure is trivial
\cite{KLM} (this would be equivalent to the triviality of the 
2-interval
inclusion for the net $\cB$).
\subsection{Normality of super\-se\-lection sectors in 
tem\-pe\-ra\-ture sta\-tes.}
\label{norm}
In this section $\gA$ will denote a local
net of von Neumann algebras on the intervals $\cI$
of $\mathbb R$ satisfying the properties $a)$ and $b)$ 
in the previous section.

With $\t$ the translation automorphism group of $\gA$, we shall say
that an endomorphism $\r$ of the quasi-local C$^*$-algebra 
$\gA$ is a {\it localized} in the 
interval $I\in\cI$ if $\r$ acts identically on $\gA(I')$, where
$I'\equiv\mathbb R\setminus I$ and, for any open set $E\in\mathbb R$,
$\gA(E)$ denotes as before the C$^*$-algebra generated by the 
$\{\gA(I): I\in\cI, I\subset E\}$.
We have also set $\cA(E)\equiv\pi_\f(\gA(E))''$.

As above, $\r$ is {\it translation covariant} if there
exists a unitary $\t$-cocycle of unitaries $u(s)\in\gA$ such that
$\Ad u(s)\cdot\t_s\cdot\r\cdot\t_{-s} = \r$.

Let $\f$ be a KMS state of $\gA$ with respect to the translation 
group. In the following $\rho$ is a translation covariant 
endomorphism of $\gA$ localized in an interval $I\in\cI$. 
By translation covariance we may assume that $I\subset (0,\infty)$. 
Our main result in this section is the following.
\begin{Thm}\label{normality} Let $\gA$ be a translation covariant 
net on $\mathbb R$ and $\f$ a KMS state of $\gA$ 
satisfying essential duality.
If $\r$ is a translation covariant localized 
endomorphism of $\gA$ with finite dimension 
$d(\rho)$, then $\r$ is normal with 
respect to $\f$, namely $\r$ extends to a normal endomorphism of
$\cM=\pi_\f(\gA)''$.

If $\f$ is an extremal KMS state, i.e. $\cM$ is a factor, then the 
extension of $\rho$ to $\cM$ has the same dimension $d(\rho)$.
\end{Thm}
Assuming that $\gA$ acts on $\cH_\f$, as above, we have:
\begin{Lemma}\label{+extension}
$\r|_{\gA(a,\infty)}$ extends to a normal endomorphism of 
$\cA(a,\infty)$ for any
$a\leq 0$.
\end{Lemma}
\begin{proof}
By translation covariance there exists a unitary $u\in\cM$
such that $\r'\equiv\Ad u \cdot \r$ is localized in an interval 
contained in $(-\infty,a)$, thus $\r=\Ad u^*\cdot\r'=\Ad u^*$ on 
$\cA(I)$ for all $I\subset(a,\infty)$, $I\in\cI$.  
It follows that $\Ad u^*$ is a normal extension of $\r$ to 
$\cA(a,\infty)$.
\end{proof}
If the endomorphism $\rho$ of $\cA$ is localized in the 
interval $I\in\cI$ 
then $\rho(\cA(I_1))\subset \cA(I_1)$
for all intervals 
$I_1$ containing $I$ 
by Haag duality. We shall say that $\r$
has {\it finite dimension} if  the index 
$[\cA(I_1):\rho(\cA(I_1))]$ is finite and independent of $I_1$
(the index is here defined for example by the Pimsner-Popa inequality 
\cite{PP}).
\begin{Lemma}\label{fi} If the endomorphism $\rho$ has finite 
dimension, then the corresponding endomorphism of $\cA(0,\infty)$ 
given by Lemma \ref{+extension} has finite dimension 
(i.e. finite index).
\end{Lemma}
\begin{proof}
Setting $\cM_n\equiv \cA(0,n)$, $n\in\mathbb N$, we have 
$\rho(\cM_n)\subset\cM_n$ for large $n$. Moreover $\cM_n$ and
$\rho(\cM_n)$ converge increasingly respectively to
$\cA(0,\infty)$ and $\rho(\cA(0,\infty))$. Thus Prop. 4
of \cite{KLM} applies and gives
\[
[\cA(0,\infty):\rho(\cA(0,\infty))]
\leq \liminf_{n\to\infty}[\cM_n:\rho(\cM_n)].
\]
\end{proof}
{\it Proof of Theorem \ref{normality}.} Let $\cB$ be the 
thermal completion of $\gA$,
thus in particular
$$
\cB(a,\infty)=\cA(\log a,\infty), \quad a>0,
$$
and denote by $U$ and $V$ the translation and dilation with respect 
to  $\cB$ as above.
We set $\a_s = \Ad V(s)$ and $\a^{(1)}_s = \Ad U(1)V(s)U(-1)$. 
By Proposition \ref{tc} then $\a^{(1)}_s$ acts on 
$\cB$ as a dilation with respect to  the point $1$, namely
\[
\a^{(1)}_s(\cB(a,\infty)) = \cB(e^s a + 1 - e^s,\infty), 
\quad a, s\in\mathbb R ,
\]
and $\a^{(1)}|_{\cB(1,\infty)}$ is the (rescaled) modular group
associated with $(\cB(1,\infty),\x)$.

As $\rho$ gives rise to a finite index  endomorphism of 
$\cB(1,\infty)$ (Lemma \ref{+extension} and \ref{fi}) there exists 
a unitary $\a^{(1)}$-cocycle $u^{(1)}(s)\in \cB(1,\infty)$
such that
\begin{equation}
\Ad u^{(1)}(s)\cdot \a^{(1)}_s \cdot\r \cdot\a^{(1)}_{-s}(X) = 
\rho(X),\quad  X\in \cB(1,\infty), \ s\in\mathbb R .
	\label{1cov}
\end{equation}
Indeed we may take $u^{(1)}$ as the Connes Radon-Nikodym cocycle
\[
u^{(1)}(s) = (D\f\cdot\Phi :D\f)_s\ ,
\]
where $\Phi$ is a normal faithful
left inverse of $\rho$ on $\cA(0,\infty)$ and
$\f$ is considered as a state on $\cA(0,\infty)$ \cite{L0}.

Now, if $\varepsilon\in(0,1)$, $\rho$ acts trivially on 
$\cB(\varepsilon,1)=\pi_\f(\cA(\log \varepsilon, 0))$ and, since 
$\f$ satisfies essential duality, the conformal thermal completion is 
strongly additive and
\[
\cB(\varepsilon,\infty)=\cB(\varepsilon,1)\vee\cB(1,\infty)\ .
\]
As $\a^{(1)}_{-s}(X)\in\cB(\varepsilon,1)$ if 
$X\in\cB(\varepsilon,1)$ and $s>0$, 
it follows that equation (\ref{1cov}) holds true for all $X 
\in\cB(\varepsilon,\infty)$,
$s>0$ .

Setting then for a fixed $s>0$
\[
\rho(X)\equiv \Ad u^{(1)}(s)\cdot \a_s^{(1)} 
\cdot\rho\cdot\a_{-s}^{(1)}(X),\quad X\in \cM \ ,
\]
this formula does not depend on the choice of $s>0$ and
provides an extension of $\rho$ to $\cM$ because of formula
(\ref{1cov}).

Now $\cB$ is a strongly additive local conformal net on $\Reali$
and (the extension of) $\r$ is a localized endomorphism of $\cB$ with finite 
dimension, hence $\r$ is M\"obius covariant \cite{GL1}, therefore the 
index of $\r(\cB(I))\subset\cB(I)$ is independent of the interval 
$I\in\cI$, provided $\r$ is localized within $I$ \cite{GL2}. This clearly 
implies the last part of the statement.
\qed
Our results then give here a version of Theorem \ref{decomposition}.
\begin{Cor}
Let $\gA$ be a net on $\Reali$ as above, $\t$ the
translation automorphism group and $\r$ a translation covariant localized 
endomorphism with finite dimension.

Then $\r$ has finite holomorphic dimension in each extremal KMS 
state $\f$ fulfilling essential duality.

The chemical potential associated with $\r$ is defined by the canonical 
splitting
\[
\b^{-1}\log d_\f(u_\r )=\b^{-1}\log d(\r) + \mu_{\r}(\f)\ ,
\]
and satisfies $\mu_{\bar\r}(\f)=-\mu_{\r}(\f)$. In particular
\[
d_{geo}(\r)=d(\r)
\]
for all irreducible $\r$, independently of $\f$.
\end{Cor}
\subsection{Extension of KMS states to the quantum double.}
\label{KMSQD}
The purpose of this subsection is to provide a description
of the chemical potential, in a low dimensional theory,
in terms of extensions of KMS states, in analogy to what described
in \cite{AHKT} in the higher dimensional case. As in our cases 
charges are not any longer associated to the dual of a compact gauge 
group, we shall replace the field algebras by a quantum double 
construction \cite{LR}, that we will perform in the C$^*$-case
in Appendix \ref{QD} for our purposes: the reader is referred to this 
appendix for the necessary notations and background.

Let then $\gA$ be a unital C$^*$-algebra with trivial centre and 
$\cT\subset{\rm End}(\gA)$ a tensor category of endomorphisms
with conjugates and sub-objects. 
Let $\a$ be a one-parameter group of automorphisms of $\gA$
and $u(\r,t)$ unitary covariance cocycle (eq. (\ref{cov})) 
which is a two-variable cocycle. We consider $\cT^{op}\subset\gA$
by setting $\r^{op}=\bar\r$ and $(\r^{op},\s^{op})=(\r,\s)^{\bullet}$. 
Let $I$ be an index set so that
$\r_i, i\in I$ is a family of inequivalent irreducible objects of 
$\cT$, one for each equivalence class. 
Then $\tilde u(\tilde\r_i,t)\equiv u(\r_i,t)\otimes 
u(\r_i,t)^{\bullet}$, where $\tilde\r=\r\otimes\r$, 
extends to a two-variable cocycle for 
$\a_t\otimes \a_t$ and $\cT\times\cT^{op}$. 
Indeed $\tilde u(\tilde\r_i,t)$ is 
independent of its choice (phase fixing) due to the anti-linearity of 
the Frobenious map $T\to T^{\bullet}$.

We may extend $\a_t\otimes \a_t$
to a one-parameter automorphism group $\tilde \a$ of $\gB$ by setting
\begin{gather}
\tilde \a_t(a)=\a_t\otimes \a_t(a), \quad a\in\tilde\gA\\
\tilde \a_t(R_i)=\tilde u(\tilde\r_i,t)^*R_i\ ,
\end{gather}
(cf. \cite{I}). If $\f$ is a KMS state for $\a_t$ on $\gA$, then 
$\f\otimes\f$ is a KMS state for 
$\a_t\otimes \a_t$ on $\tilde\gA\equiv\gA\otimes\gA$.

\begin{Prop}
Let $\f$ be an extremal KMS state for $\a_t$. 
Then $\f\otimes\f\cdot\varepsilon$ extends to 
a KMS state of $\pi_{\tilde\f}(\gB)''$ with respect to
$\tilde\a_t\cdot\theta_t$ for some one-parameter automorphism
group of $\pi_{\tilde\f}(\gB)''$ leaving $\tilde\gA$ pointwise fixed, 
if and only if each $\r\in\cT$ is 
normal with respect to $\f$. 
\end{Prop}
\begin{proof} If $\tilde\f\equiv\f\otimes\f\cdot\varepsilon$ 
is KMS state, then $\tilde\f$ is
a separating state. If $R\in\gB$ is a non-zero multiple of an 
isometry, then also 
$\psi = \tilde\f(R\cdot R^*)$ is a separating positive functional,
which is quasi-equivalent to $\tilde\f$ because $R^*\x_{\tilde\f}$
is a separating vector for $\pi_{\tilde\f}(\gB)''$.

Take $R=R_i$, the element of $\gB$ that implements $\tilde\r_i$. 
If $x\in\tilde\gA$ we have
\[
\psi(x)=\tilde\f(R_i xR_i^*)=\tilde\f(\tilde\r_i(x)E_i)
=\f\otimes\f(\varepsilon(\tilde\r_i(x)E_i)=\f\otimes\f\cdot\tilde\r_i(x)
\]
where $E_i=R_iR_i^*$, thus $\varepsilon(E_i)=1$. Thus
$\psi=\f\otimes\f\cdot\r_i$, hence
\[
\pi_{\f\otimes\f\cdot\tilde\r_i}\simeq \pi_{\psi}|_{\tilde\gA}
\simeq\pi_{\tilde\f}|_{\tilde\gA}\simeq \pi_{\f\otimes\f}\ ,
\]
and this implies that $\r_i$ is normal with respect to $\f$.

Conversely, let us assume that each $\r_i$ is normal with respect 
to $\f$ and still denote by $\tilde\r_i$ the extension of $\tilde\r_i$ to 
the von Neumann algebra $\tilde\cM\equiv\pi_{\f\otimes\f}(\tilde\gA)''$.
By rescaling the parameter we may set the inverse temperature 
$\b$ equal to $-1$.
We may then define the *-algebra $\gB_0$ as in the appendix with 
$\gA$ replaced by its weak closure 
$\tilde\cM$ (but the $C_{ij}^k$
are still defined with respect to $\tilde\gA$).

Let $\Psi_{\tilde\r_i}=
C_{\bar i i}^{0*}\bar{\tilde\r_i}(\cdot)C_{\bar i i}^0$, 
thus $\Psi_{\tilde\r_i}$ is
a (not necessarily standard) non normalized left inverse of $\tilde\r_i$.
With $V(\tilde\r_i,t)=(D\f\otimes\f\cdot\Psi_{\tilde\r_i}:D\f\otimes\f)_t$, 
define a one-parameter automorphism group
$\g$ of $\gB_0$ extending $\a\otimes\a$ as was done for $\tilde\a$, but
using $V$ instead of $U$, thus, with obvious notations, 
$\g_t(R_i)=V(\tilde\r_i,t)^*R_i$. By using
the holomorphic properties of the Connes cocycles and the corresponding 
two-variable
cocycle property (that can be checked similarly as in \cite{L4}), 
it can be seen by elementary calculations that
$\tilde\f$ is KMS with respect to $\g$. Clearly both $\g$ and $\tilde\a$
extend to $\cN=\pi_{\tilde\f}(\gB_0)''$, indeed the extension of $\g$
is the modular group with respect to $\tilde\f$, thus $\g$ and 
$\tilde\a$ commute
and $\theta_t = \tilde\a_t\cdot\g_{-t}$ is a one-parameter group of $\cN$
leaving $\tilde\cM$ pointwise fixed. 
\end{proof}
Let $G$ be the ``gauge group'', namely the group of all automorphisms 
of $\gB$ that leave $\tilde\gA$ pointwise invariant, and $Z(G)$ the 
centre of $G$. Any extension $\hat\a_t$ of $\a_t\otimes\a_t$ to a 
one parameter automorphism group of $\gB$ is clearly given by
$\tilde\a_t=\hat\a_t\cdot\theta_t$
where $\theta_t\in G$.
Reasoning as in the above proof one has the following.
\begin{Cor}
Let $\f$ and $\psi$ be extremal KMS states for $\a$. Assume that all 
$\r\in\cT$ are normal with respect to both $\f$ and $\psi$ and that if $\r$ 
is irreducible the normal extensions of $\r$ with respect to 
$\f$ and $\psi$ are 
still irreducible. Then 
$\f\otimes\psi\cdot\varepsilon$ is a KMS state of $\gB$ with respect 
to $\tilde \a\cdot\theta$, with $\theta$ a one-parameter subgroup of 
$Z(G)$.

In particular this is the case if $\f$ and $\psi$ are extremal KMS states
with respect to translations of a net of local von Neumann algebras as 
in Section \ref{CLD}, satisfying essential duality.
\end{Cor}
\noindent
The one-parameter group $\theta$ is related to 
the chemical potential associated 
with the charge $\r_i$ by
\[
\theta_t(R_i)=e^{-i\mu_{\r_i}(\f|\psi)t}R_i \ .
\]
\section{Roberts and Connes-Takesaki cohomologies.}
\label{coho}
Roberts \cite{R} has obtained a geometrical description of the theory of 
superselection sectors by considering a non-abelian local cohomology 
naturally associated with a net $\gA$. Rather than stating here his 
formal definitions, that can be obtained from 
different geometrical pictures, we will recall here the main idea 
underlying the construction of the first cohomology semiring
$H^1_R(\gA)$. Here $\gA$ is a net of von Neumann algebras on the 
Minkowski spacetime $\mathbb M$, fulfilling the properties stated in 
Section \ref{QFTM}, but the construction may be extented to more general 
globally hyperbolic spacetimes, see \cite{GLRV}.

Given a representation $\pi$ of $\gA$ localizable in all double 
cones we may, as usual, fix a double cone $\cO\in\cK$ and choose an 
endomorphism $\r_{\cO}$ of $\gA$ localized in $\cO$ and equivalent 
to $\pi$. For each $\cO_1\in\cK$ we choose a unitary intertwiner 
$u_{\cO,\cO_1}$ between $\r_{\cO}$ and $\r_{\cO_1}$ and set
\[
z_{\cO_2,\cO_1}\equiv u^*_{\cO_2,\cO}u_{\cO,\cO_1}\quad 
\cO_1,\cO_2\in\cK\ .
\]
By duality, $z_{\cO_2,\cO_1}$ belongs to $\gA^d(\cO_1\cup \cO_2)$,
in particular $z_{\cO_2,\cO_1}\in \gA(\tilde\cO)$ if $\tilde\cO$ is a 
double cone containing $\cO_1\cup \cO_2$, and satisfies the cocycle 
condition
\[
z_{\cO_3,\cO_1}=z_{\cO_3,\cO_2}z_{\cO_2,\cO_1}\ .
\]
The choice of $\r_{\cO}$ and $u_{\cO,\cO_1}$ is not unique, yet 
different choices gives cohomologous cocycles under a natural 
equivalence relation.

One can then formalize the definition of $H^1_R(\gA)$. 
The relevant point is that the map
\[
\text{Superselection sectors}\longrightarrow H^1_R(\gA)
\]
is invertible, namely each localized cocycle $z$ arises from a 
localized endomorphism $\r$ localized in a given $\cO\in\cK$, which 
is well-defined by the formula
\[\label{endo}
\r(X)=z_{\cO_1,\cO}Xz^*_{\cO_1,\cO},\quad \forall X \in 
\gA(\tilde\cO)\ ,
\]
where $\tilde\cO$ is a double cone containing $\cO$ and 
$\cO_1$ is any double cone contained in $\tilde\cO'$.

As $H^1_R(\gA)$ is in one-to-one correspondence with the 
superselection sectors, it is endowed with a ring structure, that can 
be expressed more directly. Analogusly $Z^1_R(\gA)$
is a tensor C$^*$-category.

Let's now recall the cohomology considered by Connes and Takesaki 
\cite{CT}. This concerns an automorphism group action $\t$ on a 
C$^*$-algebra $\gA$, in our case the translation automorphism group 
on the quasi-local C$^*$-algebra $\gA$. Restricting to unitary 
cocycles, a map $z:\mathbb R^4\to\gA$, taking 
values in the unitaries $\gA$, is a cocycle if satisfies the 
cocycle equation
\begin{equation}\label{CTcoc}
z(x+y)=z(x)\t_x(z(y)),\quad x,y\in\mathbb R^4\ ;
\end{equation}
two cocycles $z$ and $z'$ are cohomologous if there exists a unitary
$u\in\gA$ such that
\begin{equation}\label{CTequiv}
z'(x)=uz(x)\t_x(u^*),\quad x\in\mathbb R^4\ .
\end{equation}
A cocycle $z$ gives rise to a perturbed automorphism group
$\t^z$ defined by $\t_x^z\equiv z(x)\t_x(\cdot)z(x)^*$. We are 
interested in the case $\t^z$ is a local perturbation, in the sense 
that, if $x$ varies in a bounded set, there exists a double cone 
$\cO\in\cK$ such that
\[
\t_x^z|_{\gA(\cO')}=\t_x|_{\gA(\cO')}\ .
\]
This amounts to define a (unitary) {\it localized cocycle} as a 
unitary map $z:\mathbb R^4\to\gA$ satisfying the cocycle 
condition (\ref{CTcoc}) and the locality condition: 
\[
\exists \d>0\; \&\; \cO\in\cK\; \text{such that}\;
z(x)\in\gA(\cO),\; \forall x\in\mathbb R^4, |x|<\d\ .
\]
Here $|x|$ is the Euclidean modulus of $x$.
By iterating the cocycle equation (\ref{CTcoc}) it follows immediatly 
from the locality condition that all the $z(x)$ lives in a common 
double cone as $x$ varies in a bounded set. Indeed the following holds.
\begin{Lemma}\label{radius} 
Let $z$ be a localized cocycle. There exists $r>0$ such that
\[
z(x)\in\gA(\cO_{r+|x|}), \forall x\in\mathbb R^4\ ,
\]
where $\cO_r$ denotes the double cone with basis the ball of radius 
$r>0$, centered at $0$, in $\mathbb R^3$.
\end{Lemma}
\begin{proof} Let $\d,r>0$ be such that $z(x)\in\gA(\cO_r)$ if $|x|<\d$. 
If $(n-1)r\leq|x|<nr$, write $x= x_1+x_2+\cdots +x_n$ with $|x_i|<r$. 
By the cocycle equation
\[
z(x)=z(x_1)\t_{x_1}(z(x_2))\t_{x_1+x_2}(z(x_3))\cdots
\t_{x_1+x_2+\cdots x_{n-1}}(z(x_n))
\]
we see that $z(x)\in\cup_{i=1}^{n-1}\gA(x_1 +\cdots + x_i + 
\cO_r)\subset\gA(\cO_{(n-1)r})\subset\gA(\cO_{r+|x|})$.
\end{proof}
We denote by $Z^1_\t(\gA)$ the set of unitary localized 
cocycles, and by $H^1_\t(\gA)$ the quotient of $Z^1_\t(\gA)$
modulo the equivalence relation (\ref{CTequiv}), with the further 
restriction that the unitary $u$ is local, namely $u$
belongs to $\gA(\cO)$ for some $\cO\in\cK$.

Now a covariant localized endomorphism gives rise to a localized 
cocycle (formula (\ref{cov.formula}), hence a covariant sectors to an 
element of $H^1_\t(\gA)$. As for the Roberts cohomology, 
the converse is true.
\begin{Prop}\label{H}
There is a natural map
\[
H^1_\t(\gA) \longrightarrow 
\text{{\rm Covariant superselection sectors}}\ .
\]
\end{Prop}
\begin{proof} We shall associate a covariant localized endomorphism 
$\r$ to a given $z\in Z^1_\t(\gA)$. Similarly as in eq. (\ref{endo}), 
we set 
\[
\r(X)\equiv z(x)Xz(x)^*, \quad X\in \gA(\cO_R), \; |x|>2R\ ,
\]
where $x$ is a vector in the time-zero hyperplane and $R>r$, with $r$ 
the radius of the double cone in Lemma \ref{radius}. We have to show 
that, for a fixed $R$, the above definition is independent of the 
choice of $x$, namely $z(x)Xz(x)^*=z(x')Xz(x')^*$ if also $|x'|>2R$.
As $\cO'_{2R}$ is connected, by iterating the procedure, 
it is enough to check this if $x'=x+y$ with $|y|<\d$. Then, by local 
commutativity,
\[
z(x')Xz(x')^*=z(x)\t_x(z(y))X\t_x(z(y))^*z(x)^* = z(x)Xz(x)^*
\]
because $\t_x(z(y))$ belongs to $\gA(\cO_r + x)$ and
$\cO_r + x\subset \cO'_R$.
\end{proof}
The map given by Prop. \ref{H} is not invertible: two cocycles that 
differ by multiplication by a one dimensional character give rise to 
the same localized endomorphism; in the irreducible case this is the 
only ambiguity, that could be eliminated by considering inner 
automorphisms rather than unitary operators. Apart from this point,
$H^1_\t(\gA)$ describes faithfully the covariant 
superselection sectors. 
Note also that, if a certain strong additivity assumption holds, 
all sectors with finite dimension are covariant  \cite{GL1}.

Although the Roberts cohomology describes the 
superselection sectors, the definition of the statistical dimension 
is not manifest within $H_R^1(\gA)$. But, passing to $H_\t^1(\gA)$, we 
discover a pairing between factor KMS states $\f$
for the time evolution 
satisfying duality and $H_\t^1(\gA)$
\[
\log\langle \f, [z] \rangle \equiv 
\log d_\f(z) = \log d(\r) + \b\mu_\r(\f)
\]
hence we have the geometrical description of the holomorphic dimension 
by the diagram
\begin{equation*}
\CD
\text{Covariant sectors} @>>>H^1(\gA)
\\ @VVV	@VVV  \\  
H_\t^1(\gA)   
@>>> \mathbb R\cup\{\infty\}
\endCD
\end{equation*}
Of course, by considering the involution $^{\bullet}$ as before, 
we obtain an expression for $d(\r)$ as well.
\section{Increment of black hole entropy as an index.}
\label{blackhole}
This section is devoted to the illustration of a
physical context where the dimension is expressed by
a quantity that includes also classical geometric data of the 
underlying spacetime. The results here below have been announced in
\cite{L7}.

We shall consider the increment of
the entropy of a quantum black hole which is represented by a globally 
hyperbolic spacetime with bifurcate Killing horizon.

The case of a Rindler black hole has been studied in \cite{L4}, the 
reader can find the basic ideas and motivations in this reference. The 
main points here are the following. Firstly we perform our analysis in 
the case of more realistic black holes spacetimes, in particular we 
consider Schwarzschild black holes. 
Secondly we shall consider charges localizable on the horizon, 
obtaining in this way quantum 
numbers for the increment of the entropy of the {\it black hole itself},
rather than of the {\it outside} region as done in the Rindler spacetime.

Indeed the restriction of quantum fields to the horizon will define a 
conformal quantum field theory on $S^1$. This key point has been
discussed in \cite{L5,GLRV} with small variations. 
We shall return on this with 
further comments. Finally, we shall deal with general KMS states, besides the 
Hartle-Hawking temperature state. In this  
context, a non-zero chemical potential can appear.

The extension of the DHR analysis of the superselection structure 
\cite{DHR} to a quantum field theory
on a curved globally hyperbolic spacetime
has been given in \cite{GLRV} and we refer to this paper the 
necessary background material. We however recall here 
the construction of conformal symmetries for the observable 
algebras on the horizon of the black hole.

To be more explicit let $\cV$ be a $d+1$ dimensional 
globally hyperbolic spacetime with a bifurcate Killing
horizon. A typical example is given by the Schwarzschild-Kruskal
manifold that, by Birkoff theorem, is the only spherically symmetric
solution of the Einstein-Hilbert equation; one 
might first focus on this specific example, as the more general
case is treated similarly. We denote by $\mathfrak h_+$ and 
$\mathfrak h_-$ the two  codimension 2 submanifolds
that constitute the horizon $\mathfrak h
= \mathfrak h_+\cup \mathfrak h_-$. We assume that the horizon 
splits $\cV$ in four connected components, the future, the past and the 
``left and right wedges'' that we denote by $\cR$ and 
$\cL$ (the reader may visualize this also in the analogous case of the 
Minkowski spacetime, where the Killing flow is a one-parameter group
of pure Lorentz transformations).

Let $\k=\k(\cV)$ be the surface gravity, namely, denoting by 
$\chi$ the 
Killing vector field, the equation on $\mathfrak h$ 
\begin{equation}
	\nabla g(\chi,\chi)=-2\k\chi\ , 
	\label{k}
\end{equation}
with $g$ the metric tensor, 
defines a function $\k$ on $\mathfrak h$, that is actually constant 
on $\mathfrak h$, as can be checked 
by taking the Lie derivative of both sides of  eq. (\ref{k}) 
with respect to $\chi$ \cite{KW}.
If $\cV$ is the Schwarzschild-Kruskal manifold, then
$$
\k(\cV)=\frac{1}{4M},
$$
where $M$ is the mass of the black hole. In this case $\cR$ is 
the exterior of the Schwarzschild black hole.

Our spacetime is $\cR$ and we regard $\cV$ as a completion of $\cR$.

Let $\gA(\cO)$ be the von Neumann algebra on a Hilbert space $\cH$
of the observables localized in the 
bounded diamond $\cO\subset\cR$. We make the assumptions
of Haag duality, properly infiniteness of $\gA(\cO)$, Borchers property 
B.
The Killing flow $\L_t$ of $\cV$
gives rise to a one parameter group of automorphisms $\a$ of the
C$^*$-algebra $\gA=\gA(\cR)$ since 
$\cR\subset\cV$ is a $\L$-invariant region. Here, as before, the 
C$^*$-algebra associated with an unbounded region is the C$^*$-algebra 
generated by the von Neumann algebras associated 
with bounded diamonds contained in the region.
\subsection{The conformal structure on the black hole horizon.}
We now consider a locally normal 
$\a$-invariant state $\f$ on $\gA(\cR)$, 
that restricts to a KMS at inverse temperature 
$\b>0$ on the horizon algebra, as we will explain. This will be the case
in particular if $\f$ is a KMS state on all $\gA(\cR)$. 

For convenience, we shall assume that the net $\gA$ is already in the 
GNS representation of $\f$, hence $\f$ is 
represented by a cyclic vector $\x$. 
Denote by $\cR_a$ the wedge $\cR$
``shifted by'' $a\in\Reali$ along, say, $\mathfrak h_+$ (see \cite{GLRV}).
If $I=(a,b)$ is a bounded interval of $\Reali_+$, we set 
$$
\cC(I)=\gA(\cR_a)''\cap\gA(\cR_b)', \quad 0<a<b\ .$$ 
We denote by $\gC(I)$ ($I\subset (0,\infty)$) the C$^*$-algebra generated 
by all $\cC(a,b), b>a>0$, with $(a,b)\subset I$. We shall also 
set $\cC(I)=\gC(I)''$.

We obtain in this way a net of von Neumann algebras on the intervals of 
$(0,\infty)$, where the Killing automorphism group $\a$ acts 
covariantly by rescaled dilations. (With the due assumptions on duality, 
$\cC(I)$ turns out to be the intersection of the von Neumann algebras 
associated with all diamonds containing the ``interval $I$'' of the
horizon $\mathfrak h_+$ \cite{GLRV}.)

We can now state our assumption: 
$\f|_{\gC(0,\infty)}$ is a KMS state with respect 
to $\a$ at inverse temperature $\b>0$. 

The net $(a,b)\to\cC(e^a,e^b)$ is obviously a net on $\mathbb R$
where $\a$ is the translation automorphism group. We are therefore in 
the setting treated in Section \ref{CLD}, whose results may be now 
applied.

In particular, by using Wiesbrock theorem \cite{W}, we now show that the 
restriction of the net to the black hole horizon $\mathfrak h_+$ has 
many more symmetries than the original net.
\begin{Thm}{\rm (\cite{GLRV}).} The Hilbert space 
$\cH_{0}=\overline{\cC(I)\x}$ 
is independent of the bounded open interval $I$.

The net $\cC$ extends to a conformal net $\mathbb R$ of 
von Neumann algebras acting on $\cH_0$, where the Killing flow 
corresponds to the rescaled dilations.
\end{Thm}
\begin{proof}
Setting $\gA(a,b)\equiv\cC(e^a,e^b)$, $b>a$, we trivially obtain a local 
net on $\mathbb R$ and $\f$ is a KMS with respect to translations. 
The conformal extension of $\cC$ is then the conformal completion of $\gA$ 
given by Th. \ref{ctc}; the additivity assumption is here unnecessary 
since $\cH_{0}=\overline{\cC(I)\x}$ is independent of $I$. A detailed proof 
can be found in Prop. 4A.2 of \cite{GLRV}.
\end{proof}
This theorem says in particular that we may compactify $\Reali$ to the 
circle $S^1$, extend the definition of $\cC(I)$ for all proper intervals 
$I\subset S^1$, find a unitary positive energy representation of the 
M\"obius group $\text{PSL}(2,\Reali)$ acting covariantly on $\cC$,
so that the rescaled dilation subgroup is the Killing automorphism
group.

If $\x$ is cyclic for $\gC$ on $\cH$, namely $\cH_0 =\cH$ 
(as is true in Rindler case for a free field, see \cite{GLRV}), then
the net $\cC$ automatically satisfies Haag duality on $\mathbb R$. 
Otherwise one would pass to the dual net $\cC^d$ of $\cC$, which is
is automatically conformal and strongly additive \cite{GLW}.
The following proposition gives a condition for $\cC$ itself 
to be strongly additive.
\begin{Prop} If $\f$ is KMS on $\gA(\cR)$ and the strong additivity
for $\gA$ holds in the sense that
\[
\cC(0,1)\vee\gA(\cR_1)''=\gA(\cR)''\ ,
\]
then $\cC$ satisfies Haag duality on $\mathbb R$.
\end{Prop}                   
\begin{proof}
Note first that there is a conditional expectation $\varepsilon:
\gA(\cR_a)''\to\cC(a,\infty)$, $a\in\mathbb R$,  
given by $\varepsilon(X)\x =
EX\x$, $X\in \gA(\cR_a)''$, where $E$ is the orthogonal projection 
onto $\cH_0$.

The case $a=0$ is clearly a consequence of Takesaki's theorem since
$\cC(0,\infty)$ is globally invariant under the modular group
of $\gA(\cR)''$. Now the translations on $\gA$ (constructed as in Prop.
\ref{tc}) restrict to the 
translations on $\cC$ and commute with $\varepsilon$ 
due to the cyclicity of $\x$ for $\gC(a,\infty)$ on $\cH_0$.
Hence $\varepsilon$ maps $\gA(\cR_a)''$ onto $\cC(a,\infty)$,
$a\in\mathbb R$.

We then have 
\begin{multline}
\cC(0,1)\vee\cC(1,\infty)=\cC(0,1)\vee\varepsilon(\gA(\cR_1)'')
\\=\varepsilon(\cC(0,1)\vee\gA(\cR_1)'')=
\varepsilon(\gA(\cR)'')=\cC(0,\infty)\ ,
\end{multline}
hence $\cC$ is strongly additive, thus it satisfies Haag duality
on $\mathbb R$ because $\cC$ is conformal.
\end{proof}
In the following we assume that $\gA$ is strongly additive.
\subsection{Charges localizable on the horizon.}
We now consider an irreducible endomorphism $\r$ 
with finite dimension of $\gA(\cR)$ that is localizable in 
an interval $(a,b)$ of $\mathfrak h_+$, $a>0$, namely 
$\r$ acts trivially on $\gA(\cR_b)$ and on $\gC(0,a)$, 
thus it restricts to a localized endomorphism of $\cC$. 

This last requirement is necessary to extend $\r$ to 
a normal endomorphism of $\pi_\f(\gA(\cR))''$, with $\f$ as above or a 
different extremal KMS state. 
\smallskip

\noindent
{\it Remark.} If we assume that the net $\gA$ is defined on all $\cV$
and that $\r$ is a 
transportable localized endomorphism of $\gA(\cV)$, then $\r$ has a
normal extension to the von Neumann algebra of $\gA(\cR)''$  
exists by transportability, cf. \cite{L4,GLRV}, 
and the strong additivity assumption is unnecessary. This case may be 
treated just as the case of the Rindler spacetime \cite{L4}, the only 
difference being the possible 
appearance a non-trivial chemical potential. We
omit the detailed discussion of this context.
\smallskip

By a result in \cite{GL2}, a transportable localized
endomorphism  with finite dimension $\r$ of $\cC$ is M\"obius 
covariant. In particular we may choose the covariance cocycle 
so that it verifies the two-variable cocycle property with respect 
to the the action of the M\"obius group, and this uniquely fixes it.
Let $\s$ be another  irreducible 
endomorphism of $\gA$ localized in $(a,b)\subset\mathfrak h_+$ and
denote by $\f_\r$ and $\f_\s$ the thermal states for the 
Killing automorphism group in the 
representation $\r$. As shown in \cite{L4}, $\f_\r=\f\cdot\Phi_\r$,
where $\Phi_\r$ is the left inverse of $\r$, and similarly for $\s$.

We now show that the dimension of $\r$ and of its restriction to $\cC$
coincide. In fact the following is true.
\begin{Lemma} With the above assumptions, if $\r$ is localized
in an interval of $\mathfrak h_+$, 
then $d(\r|_{\gC(0,\infty)})$ has a normal extension to the weak closure
$\cC(0,\infty)$ with dimension $d(\r)$.
\end{Lemma}
\begin{proof} If $\r$ is localized in the interval $(a,b)$ ($b>a>0$) of
$\mathfrak h_+$, then clearly $\r$ restricts to the von Neumann algebra
$\cC(c,d)$ for all $d>b>a>c$ as it acts trivially on $\gA(\cR_c)$ and 
by duality for $\cC$. 
Hence $\r$ restricts to the C$^*$-algebra $\gC(0,\infty)$,
and then it extends to $\cC(0,\infty)$ by Theorem \ref{normality}.

Let $\bar\r$ be also localized in the interval $(a,b)$. If $R,\bar R$ are
a standard solutions for the conjugate equation of $\r$  and $\bar\r$,
then $R,\bar R\in\gA(\cR_b)'\cap\gA(\cR)''=\cC(0,b)$ and commute with
$\cC(0,a)$, hence they belong to $\cC(a,b)$ by strong 
additivity \cite{GLRV}. Conversely, if $R,\bar R$ are
a standard solution for the conjugate equation of the restriction of
$\r$ and $\bar \r$, then $R$ and $\bar R$ belongs to $\cC(a,b)$ by
the Haag duality for $\cC$, hence the conjugate equation is valid
for all elements of $\cC(0,b)\vee\gA(\cR_b)''=\gA(\cR)$. This implies 
the dimension is the same for $\r$ and its restriction to $\cC$.
\end{proof}
The increment of the free energy between the thermal equilibrium 
states $\f_\r$ and $\f_\s$ is expressed as in \cite{L4} by
$$
F(\f_\r|\f_\s)=\f_\r(H_{\r\bar{\s}})-\b^{-1} S(\f_\r|\f_\s)
$$
Here $S$ is the Araki relative entropy and 
$H_{\r\bar{\s}}$ is the Hamiltonian on $\cH_0$ 
corresponding to the composition of the 
charge $\r$ and the charge conjugate to $\s$ 
in $\cC(-\infty, 0)$ as in \cite{L4}; 
it is well defined as $e^{itH_{\r\bar{\s}}}
e^{-itH_{\iota}}$ is the two-variable cocycle.
In particular, if $\s$ is the identity 
representation, then $H_{\r\bar{\s}}=H_\r$ is the Killing
Hamiltonian in the representation $\r$.

The analysis made in Section \ref{abschem} works also in this context,
and in fact the formulae there can be written in the case 
of two different KMS states $\f_\r$ and $\f_\s$.
In particular
\begin{equation}\label{IFE}
F(\f_\s|\f_\r)=\frac{1}{2}\b^{-1}(S_c(\s)-S_c(\r)) + \mu(\f_\s|\f_\r),
\end{equation}
where $\mu(\f_\s|\f_\r)=\b^{-1}\log d_\f(u_\s )
- \b^{-1}\log d_\f(u_\r )$. Here $S_c(\r)=\log d(\r)^2$
is the conditional entropy associated with $\r$ (see \cite{PP})
and the above formula gives a canonical splitting for the incremental 
free energy.
\subsection{An index formula.}
\label{black}
We now consider the Hartle-Hawking state $\f$, see \cite{W}. 
In several cases this is the unique KMS state for the Killing 
evolution. The corresponding Hawking temperature is is related to the 
surface gravity of $\cR$:
\[
\b^{-1}= \frac{\k(\cR)}{2\pi}\ .
\]
\begin{Thm}\label{rel.free} With $\f$ the Hartle-Hawking state,
if $\r$ and $\s$ are  localizable as above, then
\[
\log d(\r)-\log d(\s)=
{\frac{\pi}{\k(\cR)}}(F(\f_\r|\f_\s)+F(\f_{\bar\r}|\f_{\bar\s})),
\]
in particular, once we fix $\r$, the exponential of
the right hand side of the above equation is proportional to an integer.
\end{Thm}
\begin{proof} 
Formula (\ref{IFE}) states that
\begin{equation}\label{ind1}
\log d(\r)-\log d(\s)=\b F(\f_\r|\f_\s) - \b\mu(\f_\s|\f_\r)\ .
\end{equation}
By the asymmetry of the chemical potential $\mu(\f_\s|\f_\r)
+\mu(\f_{\bar\s}|\f_{\bar\r})=0$, thus
\begin{equation}\label{ind2}
\log d(\r)-\log d(\s)=\b F(\f_{\bar\r}|\f_{\bar\s})+\b\mu(\f_\s|\f_\r)\ .
\end{equation}
Summing up equations (\ref{ind1}) and (\ref{ind2}) and setting 
$\b/2=\pi/\k(\cR)$ we obtain 
\[
\log d(\r)-\log d(\s)=\log d_{geo}(\r)-\log d_{geo}(\s)=
{\frac{\pi}{\k(\cR)}}(F(\f_\r|\f_\s)+F(\f_{\bar\r}|\f_{\bar\s})) .
\]
\end{proof}
We have thus expressed the analytical index $\log d(\r)-\log d(\s)$
in terms of a physical quantity, the incremental
free energy, and the quantity $\k(\cR)$ associated with the geometry
of the spacetime.
\begin{Cor}
$$
F(\f_\s|\f_\r)+F(\f_{\bar\s}|\f_{\bar\r})=\b^{-1}(S_c(\s)-S_c(\r))
$$
where $S_c(\r)$ denotes the conditional entropy of the sector $\r$.
\end{Cor}
\begin{proof} Immediate.
\end{proof}
The above corollory shows that the part of the incremental free energy 
which is independent of charge conjugation is proportional to the 
increment of the conditional entropy.
\section{On the index of the supercharge operator.}
\label{supercharge}
The purpose of this section is to make a few remarks to interpret 
the statistical dimension of a superselection sector as the Fredholm 
index of an operator associated with the supercharge operator 
in a supersymmetric theory.

Let $\gF(\cO)$ be the von Neumann algebra on a Hilbert space $\cH$ 
generated by the fields localized in the region $\cO\in\cK$ of 
a spacetime $\mathbb M$, say the Minkowski spacetime, 
in a Quantum Field Theory, as in Section \ref{QFTM}, in the vacuum 
representation. Let $\gF = \cup_{\cO\in\cK}\gF(\cO)^{-}$ be the 
quasi-local C$^*$-algebra and $\g:g\in G\to\g_g\in\text{Aut}(\gF)$ an 
action of a compact group of internal symmetries with $g_0$ an 
involutive element of the center of $G$ providing a grading 
automorphism $\g_0=\g_{g_{0}}$ with normal commutation relations:
\[
F_1 F_2 \pm F_2 F_1= 0,\quad F_i\in \gF(\cO_i),\quad\cO_1\subset\cO'_2\ ,
\]
where $F_1 , F_2\in \gF_{\pm}\equiv\{F:\g_0(F)=\pm F\}$, and the $+$ sign
occurs iff both $F_1$ and $F_2$ belong to $\gF_-$.

With $\Gamma$ the canonical selfadjoint unitary on $\cH$ implementing 
$\g_0$, the Hilbert space decomposes according to the eigenvalues of 
$\Gamma$
\begin{equation}\label{pm}
\cH=\cH_+\oplus\cH_-\ .
\end{equation}
As is known \cite{DHR}, each irreducible representation $\pi\in\hat G$ 
gives rise to a superselection sector $\r_{\pi}$ of the observable 
algebra $\gA\equiv\gF^G$: given a Hilbert space of isometries 
$H_{\pi}\subset\gF(\cO)$ carrying the representation $\pi$, one 
has a covariant endomorphism $\r_{\pi}$ of $\gA$ 
\begin{equation}\label{r1}
\r_{\pi}(X) = \sum_{i=1}^d v_i X v_i^*, \quad X\in \gF \ ,
\end{equation}
where $\{v_1,v_2,\dots v_d\}$ is an orthonormal basis of $\cH_{\pi}$
and 
\[
d=\text{dim}(\pi)=d_{DHR}(\r_{\pi})\ .
\]
The unitary representation of $G$ implementing $\g$ gives raise to a 
decomposition of $\cH$
\[
\cH=\bigoplus_{\pi\in\hat G}\cH_{\pi}
\]
and, denoting by $\pi_0$ the identity representation of $\gA$ on 
$\cH$, one has the unitary equivalence
\[
\pi_0 |_{\cH_\pi} \approx 
\text{dim}(\pi)\pi_0\cdot\r_\pi |_{\cH_{\iota}}\ .
\]
Let $H$ be the Hamiltonian of $\gF$, namely the generator of the time 
evolution on $\cH$. We now assume the existence of a supersymmetric 
structure on $\gF$, namely there exists a supercharge operator $Q$, an 
odd selfadjoint operator with $Q^2=H$; corresponding to the 
decomposition (\ref{pm}) of $\cH$, $Q$ can be written as
\begin{equation*}\begin{pmatrix} 
0 & Q_+\\ Q_- & 0 
\end{pmatrix}\end{equation*}
with $Q_- = Q_+^*$, where $Q_{\pm}:\cH_{\pm}\to\cH_{\mp}$. In case
(not assumed here) that $e^{-\b H}$ is a trace class operator, 
$\forall \b>0$, one has the well-known formula
\[
\text{Index}(Q_+)= \text{Tr}(\Gamma e^{-\b H}),\: \b>0 \ ,
\]
with $\text{Index}(Q_+)$ the Fredholm index of $Q_+$, that turns out 
to coincide with the dimension of the kernel of $H$, thus 
$\text{Index}(Q_+)=1$ if there exists a unique vacuum vector.

Now fix an irreducible localized endomorphism $\r=\r_{\pi}$ 
of $\gA$ as above and let $H_\r$ be the Hamiltonian in the 
representation $\r$, namely the generator of the unitary time 
evolution in the representation $\r$. The Hamiltonian $H_\r$ is not 
supersymmetric, in the sense that there exists no odd square root of 
$H_\r$, as such an operator would interchange Bose and Fermi sectors and 
cannot map the Hilbert space of $\r$ (which is either contained in 
$\cH_+$ or in $\cH_-$) into itself. Yet, we may define a 
supercharge operator $Q_\r$ on the global Hilbert space $\cH$ by 
setting
\[
Q_\r = \sum_i v_i Q v_i^*\ ,
\]
with the $v_i$'s forming a basis of $H_{\pi}$ as above.
\begin{Lemma} $Q_\r^2 |_{\cH_0}= H_\r$.
\end{Lemma}
\begin{proof} Indeed  
\begin{equation}\label{ss}
Q_\r^2 = (\sum_i v_i Q v_i^*)^2= \sum_i v_i Q^2 v_i^* =
\sum_i v_i H v_i^*\ .
\end{equation}
Because of the expression (\ref{r1}), one checks easily that 
$U(t)\r(X)U(-t)= \r(\a_t(X))$, $X\in\gA$, where
$U(t)=\sum_i v_i e^{itH} v_i^*$ and $\a$ is the one-parameter 
automorphism group. Since $H$ commutes with $\g_G$, it follows that
$\cH_0$ is an invariant subspace for $\sum_i v_i H v_i^*$, thus the 
latter restricts to $H_\r$ on $\cH_0$ and formula (\ref{ss}) implies
the statement.
\end{proof}
We now assume that $\r$ is a Bose sector, thus $H_{\pi}\subset \gF_+$, 
namely each $v_i$ commutes with $\Gamma$ and thus preserves 
the decomposition (\ref{pm}) of the Hilbert space so that
\begin{equation*}
Q_\r=\begin{pmatrix} 0 & Q_{\r+}\\ Q_{\r-} & 0 \end{pmatrix}
\end{equation*} 
where in particular $Q_{\r+}=\sum_i v_i Q_+ v_i^*$. 
Since the restriction of $v_i Q_+ v_i^*$ to $v_i v_i^*\cH_+$ is 
unitarily equivalent to $Q_+$, clearly $Q_{\r+}$ is unitarily 
equivalent to $Q_+\oplus\dots\oplus Q_+$ ($d$ times) and therefore
\[
\text{Index}(Q_{\r+}) = d(\r)\cdot\text{Index}(Q_+)\ .
\]
If $\s$ is another sector as above, we thus have the formula
\[
\frac{d(\r)}{d(\s)}=\frac{\text{Index}(Q_{\r+})}{\text{Index}(Q_{\s+})}
\]
showing an interpretation of the (statistical) dimension as a 
multiplicative relative Fredholm index.

It remains to provide a model where the above structure can be realized. 
To this end let $\gA_b$  and $\gA_f$ be the nets of local algebras 
associated with the free scalar Bose field and the free Fermi-Dirac field 
on the Minkowski spacetime, or on the cylinder spacetime. 
Then $\gF = \gA_b\otimes\gA_f$ is equipped with a 
supersymmetric structure (see e.g. \cite{C,H,JLO}), where the grading 
unitary is $1\otimes (-1)^N$, with $(-1)^N$ the even-odd symmetry.
We fix a positive integer $n\geq 2$ and consider
\begin{equation}\label{ntensors}
\tilde\gF\equiv\gF\otimes\gF\otimes\cdots\otimes\gF\ , (n\,\text{factors}),
\end{equation}
as our field algebra, with the gauge group $G=\mathbb P_n$ acting 
by permuting the order of the bosonic algebra $\gA_b\otimes\cdots\gA_b$ in
(\ref{ntensors}). The only thing to check is the the existence of a 
supersymmetric structure, which is given by the following lemma.
\begin{Lemma} The tensor product of two supersymmetric QFT nets is 
supersymmetric.
\end{Lemma}
\begin{proof} To simplify notations we shall consider the tensor 
product $\tilde\gF=\gF\otimes\gF$ of the same net $\gF$ by itself. 
The decomposition (\ref{pm}) of the Hilbert space $\cH$ of $\gF$ 
gives a decomposition of $\tilde\cH=\cH\otimes\cH$
\begin{align}
\tilde\cH_+ &= (\cH_+\otimes\cH_+)\oplus(\cH_-\otimes\cH_-) \\
\tilde\cH_- &= (\cH_+\otimes\cH_-)\oplus(\cH_-\otimes\cH_+) \ ,
\end{align}
and we can define the operators $\tilde 
Q_{\pm}:\tilde\cH_{\pm}\to\cH_{\mp}$
\begin{align}
\tilde Q_+ &= (Q_+\otimes 1 + 1\otimes Q_+) \oplus (Q_-\otimes 1 -
 1\otimes Q_-)\\
 \tilde Q_- &= (Q_+\otimes 1 + 1\otimes Q_-) \oplus (Q_-\otimes 1 -
 1\otimes Q_+) \ .
\end{align}
The tensor product Hamiltonian $\tilde H = (H\otimes 1)\oplus 
(1\otimes H)$ is equal to $\tilde Q_-\tilde Q_+\oplus\tilde Q_+\tilde 
Q_-$.
\end{proof}
It should be noticed that the choice of the tensor product supercharge 
$\tilde Q$ is not canonical: having chosen $\tilde Q$ we get another 
supercharge $\tilde Q_g$ by permuting the order of the tensor 
product factors with a permutation $g$.
\smallskip

\noindent
{\it Remark.} If we apply the above procedure to chiral conformal 
field theory, with $G$ a finite group of internal symmetries, 
the observable algebra $\gA=\gF^G$ has an extra family of 
sectors $\{\s_i\}_i$, beside the above ones
$\{\r_{\pi}\}_{\pi\in\hat G}$ \cite{KLM}; indeed
$\sum_{\pi}d(\r_{\pi})^2= |G|$, while
\[
\sum_{\pi\in\hat G}d(\r_{\pi})^2 +\sum_i d(\s_i)^2 = |G|^2 \ .
\]
The dimension $d(\s_i)$ is not necessarily integral and therefore 
cannot be related to a Fredholm index. Our formulae, nevertheless, 
still make sense.
\section{Outlook. Comparison with the JLO theory.}
This section contains a tentative proposal to analyze the 
superselection sectors by noncommutative geometry. Although
it is in a primitive form, we hope it will provide an insight to the
structure.
\subsection{Induction of cyclic cocycles.}
Let $\gA$ be a $\mathbb Z_2$-graded unital pre-C$^{*}$-algebra, with 
C$^*$ completion $\bar\gA$. Morphisms of $\gA$ will be assumed
to be bounded, i.e. to extend to $\bar\gA$.  
We assume that $\gA$ is a Banach $^*$-algebra with respect to a norm 
$|\!|\!|\cdot|\!|\!|$ preserved by the grading $\g$, so that 
$\gA$ is a graded Banach $^*$-algebra.

We briefly recall some basic definitions about  Connes \cite{C} entire 
cyclic cohomology. Let $\cC^n(\gA)$ be the  Banach 
space of the $(n+1)$-linear functionals of $\gA$ with finite norm
\begin{equation*}
|\!|\!|f_n|\!|\!|=\sup_{|\!|\!|a_i|\!|\!|\leq1}|f_n(a_0,a_1,\dots,a_n)|
\end{equation*}
and let $\cC(\gA)$ be the space of the entire cochains, namely 
the elements of $\cC(\gA)$ are the
sequences $f=(f_0,f_1,f_2,\dots)$, $f_n\in\cC^n(\gA)$, such that
\begin{equation}
	\|f\|=\sum_{n\geq 0}\sqrt{n!}|\!|\!|f_n|\!|\!|z^n
	\label{}
\end{equation}
is an entire function of $z$.

The grading $\g$ lifts to $\cC(\gA)$, and let $\cC_{\pm}(\gA)$ denote the 
corresponding splitting of the entire cochains.

Denoting by $\cC^e$ and $\cC^o$ the spaces of the even and odd
entire co\-chains $(f_0, f_2, f_4,\dots)$ and $(f_1, f_3, f_5,\dots)$,
the entire co\-ho\-mo\-lo\-gy groups $H^e_+(\gA)$ and $H^o_+(\gA)$  are 
the ones associated with the complex 
\begin{equation}
	\dots \to\cC^e_+\overset{\partial}{\to}\cC^o_+
	\overset{\partial}{\to}\cC^e_+\to \dots
	\label{}
\end{equation}
where the coboundary operator is $\partial= b+B$
\begin{multline*}
	(Bf)_{n-1}(a_0,a_1,\dots,a_{n-1})=
	\sum_{j=0}^{n-1}(-1)^{(n-1)j}f_n(1,a_{n-j}^{\g},
	\dots,a_{n-1}^\g,a_0\dots,a_{n-j-1})\\
	+(-1)^{n-1}f_n(a_{n-j}^\g,\dots,a_{n-j-1},1)	
\end{multline*}
 $$
 (bf)_{n+1}=\sum_{j=0}^{n}f_n(a_0,\dots,a_j a_{j+1},\dots,a_{n+1})
 +(-1)^{n+1}f_n(a_{n+1}^{\g}a_{0},a_1,\dots,a_n).
$$
If $\gB$ is another $^*$-Banach algebra, and graded 
pre-C$^*$-algebra, and $\r:\mathfrak A\to\mathfrak B$ a homomorphism 
of $\mathfrak A$ into $\mathfrak B$, any cyclic cocycle $\t$ of 
$\mathfrak B$ has a pull-back to a cyclic cocycle $(\t\cdot\r)$ of 
$\mathfrak A$
\begin{equation*}
 (\t\cdot\r)_n(a_0,a_1,\dots ,a_n)=\t_n(\r(a_0),\r(a_1),\dots,\r(a_n)).
\end{equation*}
In particular, let $\gB$ be contained in $\gA$, where both 
$\bar\gA$ and $\bar\gB$ are unital with trivial center, 
and let $\l$ be a canonical endomorphism of $\gA$ with respect to a
$\gB$, namely eq. (\ref{can2}) holds with $T\in\gA$ and $S\in\gB$, 
see the appendix. 
In this case, if $(\t_n)$ is a cyclic cocycle of
$\gB$, its pull-back to $\gA$ via $\l$ will be called the {\it 
induction} of $(\t_n)$ to $\gA$.

\begin{Prop} The induction gives rise to a well-defined map from 
$H^1(\gB)$ to $H^1(\gA)$, independently of the choice of $\l$.
\end{Prop}
\begin{proof} Immediate by the Proposition \ref{uce} and the fact that 
inner automorphisms give the identity map in cyclic cohomology 
\cite{C}.
\end{proof}
Let $\cT$ be a tensor category, with conjugates and subobjects, 
of bounded endomorphisms of $\gA$
and let $\r$ be an object of $\cT$. If $(\t_n)$ is a cyclic cocycle
of $\gA$, then we obtain a cyclic cocycle $(\t^{\r}_n)\equiv
(\t_n\cdot\r^{-1})$ on $\gB=\r(\gA)$ as pull back by $\r^{-1}$,
hence a cyclic cocycle of $\gA$ by inducting it up to $\gA$.

\begin{Prop} The induction of  $(\t_n\cdot\r^{-1})$ from $\r(\gA)$ 
to $\gA$ is equivalent to $(\t_n\cdot\bar\r)$, namely to
the pull-back of $(\t_n)$ via the conjugate charge $\bar\r$.
\end{Prop}
\begin{proof} Since $\lambda=\r\bar\r$ is the canonical endomorphism
of $\gA$ into $\r(\gA)$ (see the appendix), the result is immediate.
\end{proof}
{\it Remark.} There is a {\it product} on the localized cyclic cocycles 
(i.e. cocycles of the form $\t\cdot \r$, for some localized endomorphism 
$\r$), by setting $(\t\cdot \r)\otimes (\t\cdot \s)\equiv \t\cdot\r\s$. 
This product passes to equivalence classes and gives a well defined 
product in cohomology. The interest in this point relies on the fact 
that in the passage from the commutative to the noncommutative case the 
ring structure in cohomology is usually lost (viewed on the K-theoretical 
side, there is no noncommutative version of the tensor product of fiber 
bundles). 
\subsection{Index and super-KMS functionals.}
Let $\gA_b$, $\gA_f$ be pre-C$^*$-algebras 
with C$^*$-completions $\bar\gA_b$, 
$\bar\gA_f$, and let $\gA\equiv\gA_b\otimes\gA_f$ be
their tensor product. We assume that $\bar\gA_f$ is
$\Interi_{2}$-graded, while the grading on $\gA_b$ 
is trivial. Let $\g$ the involutive automorphism of $\gA$
giving the grading on $\gA$, trivial on $\gA_b$.
Let $\a\equiv\a^{(b)}\otimes\a^{(f)}$ be a one-parameter group of 
automorphisms of $\bar\gA$ leaving  invariant $\gA$ 
and $\d$ be an unbounded odd 
derivation of $\bar\gA$, with 
a dense $^*$-algebra $\gA_\a=\gA_{\a,b}\otimes\gA_{\a,f}\subset\gA$ 
contained in the domain
of $\d$, so that 
$\d$ is a square root of the generator $D$ of $\a$
\[
\d^2=D\equiv \frac{\textnormal{d}}{\textnormal{d}t}\a_t|_{t=0}
\quad {\rm on}\;  \gA_{\a}\ ,
\]
cf. \cite{JL,Ki}.
We consider now a super-KMS functional $\f$ at inverse 
temperature $\b=1$,
namely $\f$ is a linear functional on $\gA$ that satisfies
equation (\ref{GKMS}) for all $a,b\in\gA$ and
$$
\f(\d a)=0, \quad a \in \gA_\a.
$$
We assume that $\f$ is bounded, as is the case for quantum field 
theory on a compact space, where $\f$ is a super-Gibbs functionals. 
At the end of this section we shall add comments on the case $\f$ 
unbounded. Associated with $\f$ there is an 
entire cyclic cocycle, the JLO cocycle \cite{JLO,K,JLW}, 
on the algebra $\gA_\a$ normed with the norm $|\!|\!|a|\!|\!|\equiv
||a||+||\d a||$. For $n$ even, it is defined as
\begin{multline}\label{JLO}
\t_n(a_0,a_1,\dots,a_n)\equiv\\
(-1)^{-\frac{n}{2}}\int_{0\leq t_1\leq\dots\leq t_n\leq 1}\!\!\!\!
\f(a_0\a_{it_1}(\d a_1^\g)\a_{it_2}(\d a_2)\dots
\a_{it_n}(\d a_n^{\g^n})){\rm d}t_1{\rm d}t_2\dots{\rm d}t_n\ .
\end{multline}
There is a special class of unitary $\a$-cocycles with finite holomorphic
dimension that correspond to bounded perturbations of the dynamics
\cite{JLW}. If $q\in \gA$ is even there is a 
perturbed structure $(\gA,\a^q,\d_q)$, with 
$\d_q= \d+ [\cdot,q]$ and
$\a^q=\textnormal{Ad}u^q\cdot\a$, where $u^q$ is unitary cocycle in 
$\gA$ given by the solution of \cite{H}
\[
-i\frac{d}{dt}u^q(t)= u^q(t)\a_t(h), \quad h=\d q+q^2;
\]
the functional 
\[
\f^q(a)=\ac\f(au^q(t))
\]
is super-KMS with respect to $\a^q$ and the topological index of $u^q$ 
is given by the Chern character, the JLO cocycle $\t^q$ associated with 
$\f^q$ evaluated at the identity
\[
d_\f(u^q)=\f^q(1)\equiv\sum_{n=0}^{\infty}(-1)^n
\frac{(2n)!}{n!}\t^q_{2n}(1,1,\dots,1)\ .
\]
Of course, the second equality is trivially true as $\d 1=0$, thus
only the first term in the series may be non-zero.

We state the result on the deformation invariance if the topological 
index in \cite{JLW}:
\begin{Thm}\textnormal{(\cite{JLW}).} If $q\in \gA_-$ 
then $\f^q(1)=d_\f(u^q)=d_\f(1)=\f(1)$.\label{definv}
\end{Thm}
In the case of Th. \ref{definv}, $(\t^q_n)$  is indeed
cohomologous to $(\t_n)$.
\smallskip 

We consider now a different class of unitary cocycles. We specialize 
to the case where $\gA_b\otimes\gA_f$ is dense in the quasi-local 
C$^*$-algebra  associated with a supersymmetric quantum field theory 
(cf. Sect. \ref{supercharge}).
Let $\cT$ be a tensor category, with conjugate and subobjects, 
of localized endomorphisms contained in End$(\gA_{\a,b})$. 
We keep the same symbol for the endomorphism $\r$ of $\gA_b$ to denote 
the endomorphism $\r\otimes\iota$ of $\bar\gA$ and assume each $\r$ 
to be $\a$-covariant, in the sense that there is a $\a$-cocycle 
of unitaries $u(\r,t)\in\gA_{\a,b}$ such that is a eq. (\ref{cov}) holds.

Suppose first that $\r$ is an automorphism. Then, setting
$\d^\r\equiv \r\cdot\d\cdot\r^{-1}$, we see that
$\d^\r - \d$ acts identically
where $\r$ acts identically, that is, by Haag duality, 
$\d^\r$ is a perturbation of $\d$
by a derivation localized in a double cone. 

In general, when $\r$ is an endomorphism, $\d^\r$ 
is defined on $\r(\gA_\a)$. We shall then define
the multilinear form $\t^\r_n$ on $\r(\gA_\a)$
\begin{multline}\label{rococ}
\t^\r_n(a_0,a_1,\dots,a_n)\equiv
(-1)^{-\frac{n}{2}}\int_{\Sigma_n}\!\!\!
\f(a_0 u(is_1)\a_{is_1}(\d^\r(a_1^\g)u(is_2))
\a_{i(s_1 + s_2)}(\d^\r(a_2)u(is_3)) \\
\dots u(is_n)\a_{i(s_1+\dots +s_n)}(\d^\r(a_n^{\g^n})
u(is_{n+1})))
{\rm d}s_1\dots{\rm d}s_{n+1}
\end{multline}
where $\Sigma_n\equiv \{(s_1,\dots,s_{n+1}): s_i\geq 0,\ s_1+\dots 
+s_{n+1}=1\}$ and $u(t)\equiv u_\r(t)$. We now assume that $\f$ 
is of the form $\f_b\otimes\f_f$, which is automatically the 
case if $\f$ is factorial with a cluster property, and that $\f_b$ 
satisfies Haag duality, so that we can apply the results in Sect. 
\ref{QFTM}.

As an irreducible object $\r$ of $\cT$ extends to
irreducible endomorphisms of the weak closure on the of $\gA_{\a,b}$
in the state $\f$ (Cor. \ref{ext4dim}), then
\begin{Prop}\label{7.4}
\[
\t^\r_n = d_\f(u_\r) \t_n\cdot\r^{-1}
\]
on $\r(\gA_{\a,b})$.
\end{Prop}
\begin{proof} Setting
\[
\f_\r = \ac\f(\cdot\, u(t))\ ,
\]
we have (see Sect. \ref{factorcase} and \cite{L4}) $\f_\r = 
d_\f(u_\r)\f\cdot\Phi_\r$, thus
\[
\f_\r |_{\r(\gA)} = d_\f(u_\r)\f\cdot \r^{-1}\ ,
\]
which is a super-KMS functional with respect to the evolution
Ad$u(t)\cdot\a_t = \r\cdot\a_t \cdot\r^{-1}$ and the superderivation 
$\d^{\r}=\r\cdot\d\cdot\r^{-1}$. 

A direct verification shows that $(\t^\r_n)$ is the JLO cocycle 
on $\r(\gA_\a)$ associated with $\f_\r |_{\r(\gA_\a)}$.
\end{proof}
Thus the above expression is well defined, if
$(\t_n)$ is well defined. In order to have a cocycle on
all $\gA_{\a,b}$ we consider the induction $(\tilde\t^\r_n)$
of $(\t^\r_n)$ to $\gA_\a$.

By the previous discussion and Prop. \ref{7.4} we then have
\[
\tilde\t^\r_n = d_\f(u_\r)\t_n\cdot\bar\r\ .
\]
We summarize our discussion in the following corollary.
\begin{Cor}\label{summing} Consider the supersymmetric structure 
as above on $\gA=\gA_b\otimes\gA_f$ with $\bar\gA_b$
the quasi-local C$^*$-algebra as in Section 2.

If $\f=\f_b\otimes\f_f$ is a supersymmetric KMS functional
satisfying Haag duality and $\r$ is a translation covariant
localized endomorphism of $\gA_b$ mapping $\gA_{\a,b}$ into itself, 
then
\[
d_\f(\r) = {\tilde\t}^\r_0(1)
\left(\equiv\sum_{n=0}^{\infty}(-1)^n
\frac{(2n)!}{n!}\tilde\t^\r_{2n}(1,1,\dots,1) \right)
\]
and
\[
d_{geo}(\r) = 
\sqrt{d_\f (u_\r) d_\f (u_{\bar\r})}=
d_{DHR}(\r)\in\mathbb N \ .
\]
\end{Cor}
\smallskip

\noindent
{\it Remark.} The discussion in this section is incomplete in two 
respects. On one hand if we consider a Quantum Field Theory on a 
compact non-simply connected low dimensional spacetime, as a chiral 
conformal net on $S^1$, then the 
quasi-local C$^*$-algebra has non-trivial center with a 
non-trivial action of the superselection structure \cite{FRS};
our results then need to be extended to this context with a further 
study, possibly with a connection with the low dimensional 
QFT topology. On the other hand, when dealing 
with quantum field theory on the Minkowski spacetime, 
as in Section \ref{QFTM}, the boundedness requirement for 
the graded KMS functional $\f$ should be omitted because of 
Prop. \ref{BL}. Thus the construction of $(\tilde\t^\r_n)$, 
in its actual form, is not satisfactory except, perhaps, 
for QFT on a contractible curved spacetime. 
This difficulty vanishes if we drop the boundedness 
requirement for $\f$, but only demand the restriction
$\f_b\equiv\f |_{\gA_b}$ to be bounded. It is then natural to 
deal only with the restricion of the JLO cocycle $(\t_n)$ to the 
Bosonic algebra $\gA_{\a,b}$ . An example of such an 
unbounded super-KMS functionals seems to occur naturally
\cite{BDL}. However a general study unbounded super-KMS 
functional does not presently exist, in particular it should be
check that the JLO formula still gives
a well defined and entire cyclic cocycle $(\t_n)$.
\appendix 
\section{Appendix. Some properties of sectors in the C$^*$-case.}
For the convenience of the reader we describe here a few facts
concerning endomorphisms of C$^*$-algebras, which are natural 
counterparts of the corresponding results in the context of factors.
\subsection{The canonical endomorphism.}
Let $\gA$ be a unital  C$^*$-algebra with trivial 
centre. An endomorphism $\l$ of $\gA$ is called \emph{canonical}
(with finite index)
if there exist isometries $T\in(\iota,\l)$ and $S\in(\l,\l^2)$ such that
\begin{equation}\begin{split}
	\l(S)S=&S^2\\
	S^*\l(T)\in\Complessi\backslash\{0\},\quad &T^*S 
	\in\Complessi\backslash\{0\}\ .
	\end{split}\label{canonical}
\end{equation}
If  $\gB\subset\gA$ is a unital C$^*$-subalgebra and $\l$ is an 
endomorphism of $\gA$ with $\l(\gA)\subset\gB$ we shall say that 
$\l$ is a canonical endomorphism of $\gA$ with respect to $\gB$ (or into 
$\gB$) if there exist
intertwiners $T\in(\iota,\l)$, $S\in(\iota |_{\gB},\l|_{\gB})$
such that
\begin{equation}
	S^*\l (T)\in\Complessi\backslash\{0\},\quad T^*S 
	\in\Complessi\backslash\{0\}\ .
	\label{can2}
\end{equation}
In this case $\l$ is clearly canonical, in the sense that equation 
(\ref{canonical}) holds, since $S\in\gB$. The converse is also true.
\begin{Prop}\label{LR} If $\l$ is canonical (eq. \ref{canonical}), 
then it is canonical with respect to a
natural C$^*$-subalgebra $\gB\subset\gA$ (eq. \ref{can2}).
\end{Prop}
\noindent
The proof is the same as in the factor case \cite{L3}: one 
defines $\gB$ as the range of  $S^*\l(\cdot)S$ and  checks by eq.
(\ref{canonical}) that $\gB$ is a C$^*$-subalgebra and that  
$S^*\l(\cdot)S$ is a conditional expectation of $\gA$ onto $\gB$.
The above definitions extend to the case of pre-C$^*$-algebras, 
in this case we assume that endomorphisms are bounded and the 
intertwiners are assumed to live in the $^*$-algebras.

Note now that if $\r$ is an endomorphism of $\gA$, a conjugate $\bar\r$ 
of $\r$, if exists, is unique as sector, 
i.e. up to inner automorphisms of $\gA$, in fact a more general 
result holds in the context of tensor categories \cite{LRo}. This 
implies the following.
\begin{Prop}\label{uce} Let $\gB\subset \gA$ be unital pre-C$^*$-algebras 
with trivial centre
and $\l_1,\l_2$ canonical endomorphisms of $\gA$ into $\gB$
(namely eq. (\ref{can2}) hold). There 
exists a unitary $u\in\gB$ such that $\l_2= \textnormal{Ad}u\cdot\l_1$.
\end{Prop}
\begin{proof} The proof could be given similarly to the one given in 
\cite{L3}, Propositions 4.1 and 4.2.
However, it is easier to observe that one can consider sectors 
and conjugate sectors between different C$^*$-algebras. Then the 
canonical endomorphism $\l$, as a map of $\gA$ into $\gB$, is just 
the conjugate sector for the embedding $\iota$ of $\gB$ into $\gA$ 
and thus the uniqueness of $\l$ modulo inner automorphisms of $\gB$ 
is just a consequence of the uniqueness of the conjugate in a tensor
2-C$^*$-category, see \cite{LR}.
\end{proof}
\begin{Prop}
If $\r$ and $\bar\r$ are conjugate as above, then
$\l=\r\bar\r$ is the canonical endomorphism of $\gA$ with respect to 
the subalgebra $\r(\gA)$.
\end{Prop}
\begin{proof} Let $R,\bar R$ be the operators in the conjugate equation
for $\r,\bar\r$; setting $T=\bar R$ and $S=\r(R)$ equation (\ref{can2})
holds true.
\end{proof}
\subsection{The quantum double in the C$^*$ setting.}
\label{QD}
We consider now a construction in the C$^*$ context, that corresponds
to the one given in \cite{LR} in the context of factors, see also 
\cite{Re2,M,KLM}.

Let $\gA$ be a unital C$^*$-algebra with trivial centre and 
$\cT\subset{\rm End}(\gA)$ a tensor category of endomorphisms
with conjugates and sub-objects. We shall denote by $\gA^{op}$ the 
opposite C$^*$-algebra and by $\tilde\gA\equiv\gA\otimes\gA^{op}$
the tensor product with respect to the maximal C$^*$ tensor norm.

Let $I$ be an index set and choose $\{\r_i\}_{i\in I}$ a family all 
irreducible objects of $\cT$, one for each equivalence class, with 
$\r_0=\iota$ and $\r_{\bar i}=\bar\r_i$ and set
$\tilde\r_i\equiv\r_i\otimes\r_i^{op}$, where $\r^{op}\equiv 
j\cdot\r\cdot j$ with $j:\gA\to\gA^{op}$ the canonical anti-linear 
isomorphism.

Let $\gB_0$ be the linear space of functions $X: I\to\tilde\gA$ with 
finite support. We consider the following product and $^*$-operation 
on $\gB_0$:
\begin{gather}
X\star Y(k)\equiv\sum_{i,j}X(i)\tilde\r_i(Y(j))C_{ij}^k\ ,\\
X^*(k)\equiv C_{k\bar k}^{0*}\tilde\r_k(X(\bar k)^*)\ .
\end{gather}
Here $C_{ij}^{k}\in(\tilde\r_k,\tilde\r_i\tilde\r_j)$ is the canonical 
intertwiner $C_{ij}^{k}=
\sqrt{\frac{d(\r_i)d(\r_j)}{d(\r_k)}}
\sum_{\ell}v_{\ell}\otimes j(v_{\ell})$ 
with $\{v_{\ell}\}_{\ell}$ any orthonormal bases in $(\r_k,\r_i\r_j)$.
Clearly $\tilde\gA$ can be identified with the subalgebra of $\gB_0$ 
of functions with support in $0$.

Setting $R_i(k)\equiv \delta_{ik}$, the $R_i\in\tilde\gA$ satisfy the 
relations
\begin{equation}\label{r}
\left\{
\begin{array}{l} R_i X=\tilde\rho_i (X)R_i,
        \qquad X\in \tilde\gA \  , \\
     R_i^*R_i =d(\rho_i)^2\  ,\\
 R_iR_j=\sum_k C^k_{ij}R_k \  , \\
    R_i^*=C^{0*}_{{\bar i}i} R_{\bar i} \ ,
      \end{array}   \right.
\end{equation}
where we have omitted the $\star$ in the product. $\gB_0$ is the 
a $^*$-algebra and every $X\in\gB_0$ as a unique expansion
\[
X=\sum_i X(i)R_i\ 
\]
where the coefficients are uniquely determined by 
$X(i)=\varepsilon (XR_i^*)$. Here $\varepsilon: 
\gB_0\to\tilde\gA$ is the conditional expectation given by 
$\varepsilon(X)=X(0)$, that can be shown to be faithful as in 
\cite{KLM}, App. A.

Extending states from $\tilde\gA$ to $\gB_0$ via $\varepsilon$ 
we obtain a faithful family of states, hence $\gB_0$ has a maximal C$^*$-norm, the 
completion under which is a C$^*$-algebra $\gB$. Now
\[
||X^*X||\equiv\sup_{\psi}\psi(X^*X)\geq
\sup_{\f_1\otimes\f_2}\f_1\otimes\f_2\cdot\varepsilon(X^*X)
=||\varepsilon(X^*X)||
\]
where $\psi$ ranges over the states of $\gB_0$ and 
$\f_1\otimes\f_2$ over
the product states of $\tilde\gA$, namely $\varepsilon$ is bounded and 
extends to a conditional expectation $\varepsilon: \gB\to\tilde\gA$.

To each $X\in\gB$ we may associate the formal expansion
\[
X=\sum_i X(i)R_i\ 
\]
where $X(i)\equiv\varepsilon (XR_i^*)$ and one has $X=0\Leftrightarrow 
X(i)=0\ \forall i\in I$.

We have $\tilde\gA'\cap\gB=\mathbb C$. Indeed if 
$X\in\tilde\gA'\cap\gB$, then for every $a\in\tilde\gA$
\[
\sum_i aX(i)R_i=aX=Xa=\sum_i X(i)R_i a= \sum_i X(i)\tilde\r_i(a)R_i\ ,
\]
thus $X(i)\in(\iota,\tilde\r_i)$, thus it follows by the irreducibility 
of $\tilde\r_i$ that $X=X(0)\in\gA'\cap\gA=\mathbb C$.

The contact with the construction in \cite{LR} is visible by
the following proposition.

\begin{Prop} If $\cT$ is rational (namely $I$ is a finite set), then
\[
\l=\bigoplus_{i\in I}\r_i\otimes\r_i^{op}
\]
is the restriction to $\tilde\gA$ of the canonical endomorphism of 
$\gB$ to $\tilde\gA$.
\end{Prop}

Although the chemical potential can be described via extension of KMS
state from $\tilde\gA$ to $\gB$, in Section \ref{KMSQD} we prefer 
to deal with a slight variation 
of the above construction. 

Indeed the definition of $\gB$ can be modified by replacing $\tilde\gA$ with
$\gA\otimes\gC$, where $\gC$ is any C$^*$-algebra with trivial centre
where $\cT$ acts faithfully, namely $\cT\subset{\rm End}(\gC)$.
Accordingly one puts in this case 
$\tilde\r_i\equiv \r_i\otimes\bar\r_i$ and $C_{ij}^k=\sum_{\ell}
\sqrt{\frac{d(\r_i)d(\r_j)}{d(\r_k)}}v_{\ell}\otimes v_{\ell}^{\bullet}$.
In particular we may take $\gC=\gA$ and indeed we specialize to 
this case in section \ref{KMSQD}.  It is  rather obvious  how 
to formulate in the above setting  the structure and the results 
obtained for the quantum double case.
\bigskip

\noindent
{\bf Acknowledgments.} Among others, we wish to thank in particular
S. Doplicher for early motivational comments, and A. Connes, K. 
Fredenhagen and A. Jaffe for inspiring conversations 
and invitations respectively at the IHES, Hamburg University 
and Harvard University, while this work was at different stages.

\end{document}